%


\input xy
\xyoption{all}




\font\sc=cmcsc10 \rm


\newcount\secnb
\newcount\subnb
\newcount\parnb
\newcount\itemnb
\secnb=0
\subnb=0
\parnb=0
\itemnb=96
\newtoks\secref
\secref={}
\newtoks\subref
\subref={}

\def\smallskip{\par\vskip 2mm}
\def\medskip{\par\vskip 5mm}
\def\goodbreak{\penalty -100}


\def\section#1{\global\advance\secnb by 1
\secref=\expandafter{\the\secnb}
\subref={}
\subnb=0\parnb=0\itemnb=96
\medskip\goodbreak
\centerline{\bf\S\the\secref.\ #1}
\smallskip\nobreak}

\def\subsection#1{\global\advance\subnb by 1
\secref=\expandafter{\the\secnb.}
\subref=\expandafter{\the\subnb}
\parnb=0\itemnb=96
\smallskip\goodbreak
\leftline{\bf\the\secref\the\subref.\ #1}
\par\nobreak}

\def\references{\medskip\goodbreak
\centerline{\bf References}
\smallskip\nobreak}


\def\tit{\global\advance\parnb by 1
\itemnb=96
\smallskip\goodbreak
\noindent\the\secref\the\subref.\the\parnb. }

\def\endtit{\par\nobreak}


\def\th{\tit{\sc Theorem}}

\def\prop{\tit{\sc Proposition}}


\def\lem{\tit{\sc Lemma}}

\def\proof{\smallskip\goodbreak\noindent{\it Proof:}\ }

\def\endproof{\hbox{\unskip\kern 6pt\penalty 500
\raise -1pt \hbox{\vrule\vbox to 6pt{\hrule width 4pt
\vfill\hrule}\vrule}}}



\def\ifundefined#1{\expandafter\ifx\csname #1\endcsname\relax }


\newtoks\refid

\def\label#1{\relax}

\def\ref#1{\ifundefined{crossref#1}\refid={??}\else
\refid={\csname crossref#1\endcsname }\fi{\the\refid}}

\def\eqnnb#1{\global\advance\itemnb by 1\relax
\writeref{crossref#1}{(\the\secref\the\subref\the\parnb.\char\the\itemnb)}\relax
\hbox{(\char\the\itemnb)}}

\def\fig#1{\global\advance\itemnb by 1\relax
\writeref{crossref#1}{(\the\secref\the\subref\the\parnb.\char\the\itemnb)}\relax
\centerline{\it Figure (\the\secref\the\subref\the\parnb.\char\the\itemnb)}}


\newcount\bibrefnb
\bibrefnb=0

\def\biblabel#1{\global\advance\bibrefnb by 1\relax
\item{\the\bibrefnb.}}

\def\refto#1{\setbox1=\hbox{\sc #1}\copy1}
\def\reftosame{\hbox to \the\wd1{\hrulefill}}

\def\bibref#1{\ifundefined{bibref#1}\refid={??}\else\refid={\csname bibref#1\endcsname}\fi\relax
\the\refid}

\def\bye{\medskip{\it Address:} Benoit Fresse,
Laboratoire J.A. Dieudonn\'e,
Universit\'e de Nice,
Parc Valrose,
F-06108 Nice Cedex 02 (France).

{\it E-mail:} fresse@math.unice.fr\vfill\supereject\end}


\def\N{{\bf N}}
\def\Z{{\bf Z}}
\def\F{{\bf F}}
\def\W{{\cal W}}
\def\E{{\cal E}}
\def\B{{\cal B}}

\def\Map{\mathop{\rm Map}\nolimits}
\def\Hom{\mathop{\rm Hom}\nolimits}
\def\rto#1#2{\smash{\mathop{\hbox to 6 mm{\rightarrowfill}}\limits^{\scriptstyle #1}_{\scriptstyle #2}}}
\def\rdouble#1#2{\raise -2pt\vbox{\baselineskip=4pt
\hbox{$\displaystyle\rto{#1}{}$}
\hbox{$\displaystyle\rto{}{#2}$}}}
\def\lto#1#2{\smash{\mathop{\hbox to 6 mm{\leftarrowfill}}\limits^{\scriptstyle #1}_{\scriptstyle #2}}}
\def\id{{id}}
\def\Spaces{{\cal S}}
\def\Alg{\mathop{{\cal A}lg}\nolimits}

\def\Sq{\mathop{{\rm Sq}}\nolimits}
\def\pt{{pt}}
\def\colim{\mathop{{\rm colim}}}
\def\Tor{\mathop{{\rm Tor}}\nolimits}
\def\op{{op}}
\def\Ho{\mathop{\rm Ho}}

\def\sgn{\mathop{\rm sgn}}


\magnification=1200

\centerline{\bf Derived division functors and mapping spaces}
\smallskip
\centerline{Benoit Fresse}
\centerline{2/1/2003}
\medskip

\medskip\centerline{\bf Introduction}

\medskip
The normalized cochain complex of a simplicial set $N^*(Y)$
is endowed with the structure of an $E_\infty$ algebra.
More specifically,
we prove in a previous article
that $N^*(Y)$ is an algebra over the Barratt-Eccles operad
({\it cf}. [\bibref{BF}]).
According to M. Mandell,
under reasonable completeness assumptions,
this algebra structure determines the homotopy type of $Y$
({\it cf}. [\bibref{Mnd}]).
In this article,
we construct a model of the mapping space $\Map(X,Y)$.
For that purpose,
we extend the formalism of Lannes' $T$ functor ({\it cf}. [\bibref{L}])
in the framework of $E_\infty$ algebras.
Precisely,
in the category of algebras over the Barratt-Eccles operad,
we have a division functor $-\oslash N_*(X)$
which is left adjoint to the functor $\Hom_\F(N_*(X),-)$.
We prove that the associated left derived functor $-\mathop{\oslash}^L N_*(X)$
is endowed with a quasi-isomorphism
$N^*(Y)\mathop{\oslash}^L N_*(X)\,\rto{\sim}{}\,N^*\Map(X,Y)$.

\medskip\centerline{\bf Summary}

\medskip{\S 1. {\it Results}}. ---
We give a detailed introduction to our results in this section.

{\S 2. {\it Cofibrant resolutions and division functors}}. ---
We recall the construction of cofibrant objects in the context of $E_\infty$ algebras.
This notion occurs in the definition of a left derived functor.

{\S 3. {\it The example of loop spaces}}. ---
We consider the case of a division functor $-\oslash K$ where $K$ is the
reduced chain complex of the circle $K = \tilde{N}(S^1)$.
We make explicit the associated left derived functor $N^*(Y)\mathop{\oslash}^L \tilde{N}_*(S^1)$
which provides a model of loop spaces $\Omega Y$.

{\S 4. {\it The case of Eilenberg-MacLane spaces}}. ---
We prove our main result in section 4.
First, we construct a quasi-isomorphism $N^*(Y)\mathop{\oslash}^L N_*(X)\,\rightarrow\,N^*\Map(X,Y)$
in the case of an Eilenberg-MacLane space $Y = K(\Z/2,n)$.
Then, we proceed by induction on the Postnikov tower of $Y$.

{\S 5. {\it The Eilenberg-Zilber equivalence}}. ---
We observe that the morphism $N^*(X\times Y)\,\rightarrow\,N^*(X)\otimes N^*(Y)$
provided by the classical Eilenberg-Zilber equivalence
is not compatible with $E_\infty$ algebra structures.
We show how to overcome this difficulty in section 5.
The results of this section
allow to define a morphism
$N^*(Y)\mathop{\oslash}^L N_*(X)\,\rightarrow\,N^*\Map(X,Y)$
functorial in $X$ and $Y$.
The functoriality property is crucial in the proof of our main theorem.

\section{Results}

\tit{\it Conventions}\endtit
We fix a ground field $\F$.
We work in the category of differential graded modules over $\F$ (for short, dg-modules).
To be more precise,
we consider either upper graded modules $V = V^*$
or lower graded modules $V = V_*$.
The relation $V_d = V^{-d}$ makes a lower grading equivalent to an upper grading.
In general, the differential of a dg-module is denoted by $\delta: V\,\rightarrow\,V$.
The cohomology of an upper graded dg-module is denoted by $H^*(V)$.

We consider the classical tensor product of dg-modules
which equips the category of dg-modules with the structure of a closed symmetric monoidal category.
According to the classical rule,
we assume that the symmetry isomorphism $V\otimes W\,\rightarrow\,W\otimes V$ involves a sign,
which we denote by $\pm$ without more specification.
The dg-module of homogeneous morphisms, denoted by $\Hom_\F(V,W)$,
is the internal hom in the category of dg-modules.
To be explicit,
a homogeneous morphism of upper degree $d$ is a map $f: V\,\rightarrow\,W$
which raises upper degrees by $d$.
The differential of $f\in\Hom_\F(V,W)$
is given by the graded commutator
$\delta(f) = f \delta - \pm\delta f$.
The dual dg-module of a dg-module, denoted by $V^{\vee}$,
is defined by the relation $V^{\vee} = \Hom_\F(V,\F)$.

\tit{\it The Barratt-Eccles operad}\endtit\label{defnOperad}
We consider the differential graded Barratt-Eccles operad,
which we denote by the letter $\E$.
The classical Barratt-Eccles operad,
which we denote by $\W$,
is the simplicial operad
whose term $\W(r)$
is the homogeneous bar construction
of the symmetric group $\Sigma_r$ ({\it cf}. [\bibref{BE}], [\bibref{BF}]).
The term $\E(r)$ of the associated differential graded operad is just the normalized chain complex of $\W(r)$.
In particular, the dg-module $\E(r)$ is a non-negatively lower graded chain complex.

To be explicit, an $n$-dimensional simplex in $\W(r)$ is an $n+1$-tuple $(w_0,\ldots,w_n)$,
where $w_i\in\Sigma_r$, for $i = 0,\ldots,n$.
We have $d_i(w_0,\ldots,w_n) = (w_0,\ldots,\widehat{w_i},\ldots,w_n)$
and $s_j(w_0,\ldots,w_n) = (w_0,\ldots,w_j,w_j,\ldots,w_n)$.
The set $\W(r)$ is equipped with the diagonal action of $\Sigma_r$.

Let us consider the case $r = 2$.
The chain complex $\E(2)_*$ is the standard free resolution of the trivial representation of $\Sigma_2$.
The elements of $\Sigma_2$
(the identity permutation and the transposition)
are denoted by $\id\in\Sigma_2$ and $\tau\in\Sigma_2$
respectively.
Let $\theta_d\in\W(2)_d$ be the simplex such that $\theta_d = (\id,\tau,\id,\tau,\id,\ldots)$.
The elements $\theta_d = (\id,\tau,\id,\tau,\id,\ldots)$
and $\tau\cdot\theta_d = (\tau,\id,\tau,\id,\tau,\ldots)$
are the non-degenerate simplicies of $\W(2)_d$.
Therefore, we have $\E(2)_d = \F_2[\Sigma_2]\,\theta_d$
where $\delta(\theta_d) = \tau\cdot\theta_{d-1} + (-1)^{d-1}\theta_{d-1}$.

The operad structure is specified by a composition product
$$\E(r)\otimes\E(s_1)\otimes\cdots\otimes\E(s_r)\,\rightarrow\,\E(s_1+\cdots+s_r)$$
which arises from an explicit substitution process of permutations
(for more details, we refer to [\bibref{BF}]).
The composite of $\rho\in\E(r)$ with $\sigma_1\in\E(s_1),\ldots,\sigma_r\in\E(s_r)$
is denoted by $\rho(\sigma_1,\ldots,\sigma_r)\in\E(s_1+\cdots+s_r)$.

The set $\W(1)$ is reduced to a point.
The associated basis element is denoted by $1\in\E(1)$,
because this element is a unit for the operad composition product.
The set $\W(0)$ is also reduced to a point.
The associated basis element is denoted by $\eta\in\E(0)$,
because this element determines unit morphisms
in the context of algebras over the Barratt-Eccles operad.

\tit{\it Algebras}\endtit
We consider algebras and coalgebras over the Barratt-Eccles operad.
As a reminder, an $\E$-algebra is a dg-module $A$ equipped with an evaluation product
$\E(r)\otimes A^{\otimes r}\,\rightarrow\,A$.
Equivalently,
an element $\rho\in\E(r)$ determines a multilinear operation
$\rho: A^{\otimes r}\,\rightarrow\,A$,
which maps an $r$-tuple of elements $a_1,\ldots,a_r\in A$
to an element denoted by $\rho(a_1,\ldots,a_r)\in A$.
The operad unit $1\in\E(1)$ is supposed to give the identity operation $1(a_1) = a_1$.
Furthermore, the composite of an operation $\rho: A^{\otimes r}\,\rightarrow\,A$
with $\sigma_1: A^{\otimes s_1}\,\rightarrow\,A,\ldots,\sigma_r: A^{\otimes s_r}\,\rightarrow\,A$
has to agree with $\rho(\sigma_1,\ldots,\sigma_r): A^{\otimes s_1+\cdots+s_r}\,\rightarrow\,A$,
the operation supplied by the composition product of the operad
(as indicated by our notation).
The evaluation product is a morphism of dg-modules.
Therefore, we have the derivation relation
$\delta(\rho(a_1,\ldots,a_r)) = (\delta\rho)(a_1,\ldots,a_r)
+ \sum_{i=1}^r\pm\rho(a_1,\ldots,\delta a_i,\ldots,a_r)$.
The evaluation product is also assumed to be invariant under the action of the symmetric group.
Therefore, we have the relation
$(w\rho)(a_1,\ldots,a_r) = \pm\rho(a_{w(1)},\ldots,a_{w(r)})$.

The basis element $\eta\in\E(0)$ yields a morphism $\eta: \F\,\rightarrow\,A$.
(This morphism is also determined by an element $1\in A$ denoted as the unit of the algebra $A$.)
In fact, the ground field $\F$ has a unique algebra structure
and is an initial object of the category of $\E$-algebras.

Dually, an $\E$-coalgebra is a dg-module $K$ equipped with a coproduct $\E(r)\otimes K\,\rightarrow\,K^{\otimes r}$.
Hence, an element $\rho\in\E(r)$ determines a co-operation $\rho^*: K\,\rightarrow\,K^{\otimes r}$.
Furthermore, the ground field $\F$ is equipped with an $\E$-coalgebra structure.
In fact, it represents the final object in the category of $\E$-coalgebras,
because the element $\eta\in\E(0)$ determines an augmentation morphism $\eta^*: K\,\rightarrow\,\F$,
for all $\E$-coalgebras $K$.

As for classical algebras, the linear dual of an $\E$-coalgebra is an $\E$-algebra.
The converse assertion holds if we assume that an $\E$-algebra is finite or is equipped with a profinite topology.

\tit{\it Cochain algebras}\endtit
The normalized chain complex of a simplicial set $N_*(X)$
is equipped with the structure of a coalgebra
over the differential graded Barratt-Eccles operad $\E$
({\it cf}. [\bibref{BF}]).

Precisely, for $X$ a simplicial set, we have a natural coproduct
$\E(r)\otimes N_*(X)\,\rightarrow\,N_*(X)^{\otimes r}$
such that the co-operation $\theta_0{}^*: N_*(X)\,\rightarrow\,N_*(X)^{\otimes 2}$
determined by the element $\theta_0\in\E(2)_0$
is the classical Alexander-Whitney diagonal.
Thus, the dual operation $\theta_0: N^*(X)\otimes N^*(X)\,\rightarrow\,N^*(X)$
agrees with the classical cup-product.

Furthermore,
because we have the derivation relation $\delta(\theta_i) = \tau\cdot\theta_{i-1} + (-1)^{i-1}\theta_{i-1}$,
the higher operation $\theta_i: N^*(X)\otimes N^*(X)\,\rightarrow\,N^{*+i}(X)$
associated to the element $\theta_i\in\E(2)$
is a representative of the classical cup-$i$-product
({\it cf}. [\bibref{Stn}]).
Consequently, if $x\in N^n(X)$ is a representative of a class $c\in H^n(X)$,
then $\theta_{n-i}(x,x)\in N^{n+i}(X)$
is a representative of the $i$th reduced square $\Sq^i(c)\in H^{n+i}(X)$.

The coalgebra augmentation $\eta^*: N_*(X)\,\rightarrow\,\F$
is equivalent to the canonical augmentation $N_*(X)\,\rightarrow\,N_*(\pt)$,
because the module $N_*(\pt) = \F$ is identified with the final coalgebra.

\tit{\it The tensor structure}\endtit\label{TensorAlgebra}
The chain complex $\E(r)$ is equipped with a diagonal $\Delta: \E(r)\,\rightarrow\,\E(r)\otimes\E(r)$.
We have explicitly $\Delta(w_0,\ldots,w_n) = \sum_{k=0}^n (w_0,\ldots,w_k)\otimes(w_k,\ldots,w_n)$.
As an example, if we consider the element $\theta_d\in\E(2)_d$ introduced in paragraph \ref{defnOperad},
then we obtain $\Delta(\theta_d) = \sum_{k=0}^d \theta_k\otimes\tau^k\cdot\theta_{d-k}$.

The diagonal of $\E$ commutes with the operad composition product of $\E$.
Consequently, if $A$ and $B$ are $\E$-algebras,
then the tensor product $A\otimes B$ is equipped with the structure of an $\E$-algebra.
The elements of $\E(r)$ operate on $A\otimes B$ through the diagonal of $\E(r)$.
To be more explicit, the operation
$\rho: (A\otimes B)^{\otimes r}\,\rightarrow\,(A\otimes B)$
associated to an element $\rho\in\E(r)$
statisfies the formula
$$\rho(a_1\otimes b_1,\ldots,a_r\otimes b_r)
= \sum\nolimits_i\pm\rho_{(1)}^i(a_1,\ldots,a_r)\otimes\rho_{(2)}^i(b_1,\ldots,b_r),$$
where $\Delta(\rho) = \sum_i\rho_{(1)}^i\otimes\rho_{(2)}^i$.

Similarly, if $A$ is an $\E$-algebra and $K$ is an $\E$-coalgebra,
then the dg-module of homogeneous morphisms $\Hom_\F(K,A)$
is equipped with the structure of an $\E$-algebra.
In fact, if $K$ is a finite dimensional module,
then $K$ is equivalent to the linear dual of an $\E$-algebra $B = K^{\vee}$.
The dg-module $\Hom_\F(K,A)$ is equivalent to the tensor product of $A$
with an $\E$-algebra $A\otimes B = A\otimes K^{\vee}$.
Therefore,
in the finite dimensional case,
the dg-module $\Hom_\F(K,A)$ is equipped with the structure of an $\E$-algebra
according to the construction of the paragraph above.

Suppose given $u_1,\ldots,u_r\in\Hom_\F(K,A)$.
In general,
for $\rho\in\E(r)$,
the homogeneous morphism $\rho(u_1,\ldots,u_r)\in\Hom_\F(K,A)$
is the sum
$$\rho(u_1,\ldots,u_r) = \sum\nolimits_i\pm\rho_{(1)}^i\cdot u_1\otimes\cdots\otimes u_r\cdot\rho_{(2)}^i{}^*$$
of the composites
$$\displaystyle K
\,\mathop{\hbox to 16mm{\rightarrowfill}}^{\rho_{(2)}^i{}^*}\,K\otimes\cdots\otimes K
\,\mathop{\hbox to 16mm{\rightarrowfill}}^{u_1\otimes\cdots\otimes u_r}\,A\otimes\cdots\otimes A
\,\mathop{\hbox to 16mm{\rightarrowfill}}^{\rho_{(1)}^i}\,A,$$
where $\Delta(\rho) = \sum_i\rho_{(1)}^i\otimes\rho_{(2)}^i$.

\tit{\sc Proposition}\endtit
{\it Let $K$ be an $\E$-coalgebra. The functor $\Hom_\F(K,-)$ has a left adjoint.
More explicitly, for any $\E$-algebra $A$, there is an $\E$-algebra $A\oslash K$
such that $\Hom_\E(A\oslash K,-)\,\simeq\,\Hom_\E(A,\Hom_\F(K,-))$.}

\smallskip
This result follows readily from the special adjoint functor theorem ({\it cf}. [\bibref{MLCat}, Chapter V]).
The functor $\Hom_\F(K,-)$ preserves all limits of $\E$-algebras
because the forgetful functor from the category of $\E$-algebras to the category of $\F$-modules
creates limits.

\tit{\it A closed model category structure}\endtit
The category of $\E$-algebras is equipped with the structure of a closed model category.
By definition,
a weak-equivalence is a morphism of $\E$-algebras $f: A\,\rightarrow\,B$
which induces an isomorphism in cohomology $H^*(f): H^*(A)\,\rightarrow\,H^*(B)$
(namely, we assume that $f: A\,\rightarrow\,B$ is a quasi-isomorphism).
A fibration is a morphism of $\E$-algebras
which is surjective in all degrees.
A cofibration is a morphism of $\E$-algebras
which has the left-lifting property with respect to acyclic fibrations.
We prove that these definitions provide a full model structure
on the category of $\E$-algebras
in [\bibref{BF}].

We remind the reader that an acyclic fibration (respectively, an acyclic cofibration)
denotes a morphism which is both a fibration (respectively, a cofibration)
and a weak equivalence.
An $\E$-algebra $F$ is cofibrant if the initial morphism $\eta: \F\,\rightarrow\,F$ is a cofibration.
Similarly, an $\E$-algebra $A$ is fibrant if the final morphism $A\,\rightarrow\,0$ is a fibration.
But, this property holds for any $\E$-algebra,
because the final object in the category of $\E$-algebras is the $0$ module.
A cofibrant resolution of an $\E$-algebra $A$
is a cofibrant algebra $F$
equipped with an acyclic fibration $F\,\rto{\sim}{}\,A$.

\tit{\it The homotopy category of $\E$-algebras}\endtit
The category of $\E$-algebras is denoted by $\E\Alg$.
The associated homotopy category is denoted by $\Ho(\E\Alg)$.
The set of algebra morphisms from $A\in\E\Alg$ to $B\in\E\Alg$ is denoted by $\Hom_\E(A,B)$.
The morphism set in the homotopy category is denoted by $[A,B]_\E$.
In general,
we have $[A,B]_\E = \Hom_\E(F,B)/\sim$,
where $F$ is a cofibrant resolution of $A$,
and where $\sim$ refers to the homotopy relation 
in the category of $\E$-algebras.
In our context, it is possible to make this relation explicit.

Precisely,
the construction of paragraph \ref{TensorAlgebra} makes the definition of a path-object
in the category of $\E$-algebras easy.
Explicitly,
for any $\E$-algebra $A$,
the module $\Hom_\F(N_*(\Delta^1),A)$
is a path object of $A$,
because we have a diagram
$$\underbrace{\Hom_\F(N_*(\pt),A)}_{=A}\,\rto{\sim}{s_0}\,\Hom_\F(N_*(\Delta^1),A)
\,\rdouble{d_0}{d_1}\,\underbrace{\Hom_\F(N_*(\pt),A)}_{=A}$$
deduced from the cosimplicial structure of $\Delta^\bullet$.
We call the $\E$-algebra $\Hom_\F(N_*(\Delta^1),A)$ a cylinder object of $A$,
because we consider tacitly the opposite category of $\E$-algebras.

Suppose given a pair of parallel morphisms $f: F\,\rightarrow\,B$ and $g: F\,\rightarrow\,B$,
where $F$ is a cofibrant $\E$-algebra.
These morphisms are homotopic $f\sim g$
if and only if there is a morphism $h: F\,\rightarrow\,\Hom_\F(N_*(\Delta^1),B)$
such that $d_0 h = f$ and $d_1 h = g$.

\tit{\it Derived functors}\label{DerivedFunctors}\endtit
Clearly, the functor $\Hom_\F(K,-): \E\Alg\,\rightarrow\,\E\Alg$ carries fibrations to fibrations
and preserves all weak-equivalences of $\E$-algebras.
By adjunction,
the functor $-\oslash K: \E\Alg\,\rightarrow\,\E\Alg$
carries (acyclic) cofibrations to (acyclic) cofibrations.
As a consequence,
the functors
$-\oslash K: \E\Alg\,\rightarrow\,\E\Alg$ and $\Hom_\F(K,-): \E\Alg\,\rightarrow\,\E\Alg$
determine a pair of derived adjoint functors
$$-\mathop{\oslash}^L K: \Ho(\E\Alg)\,\rightarrow\,\Ho(\E\Alg)
\qquad\hbox{and}
\qquad \Hom_\F(K,-): \Ho(\E\Alg)\,\rightarrow\,\Ho(\E\Alg).$$

In this article, the notation $A\mathop{\oslash}^L K$ refers also to a representative of the image of $A\in\E\Alg$
under the left derived functor of $-\oslash K: \E\Alg\,\rightarrow\,\E\Alg$.
Such a representative is provided by an $\E$-algebra $F\oslash K\in\E\Alg$,
where $F$ is a cofibrant resolution of $A$.
According to this convention,
the right derived functor of $\Hom_\F(K,-): \E\Alg\,\rightarrow\,\E\Alg$
does not differ from itself,
because all $\E$-algebras are fibrant objects.

\th\label{resolutionMapSpace}\endtit
{\it Let $X$ and $Y$ be simplicial sets.
We assume that $\pi_n(Y)$ is a finite $p$-group for $n\geq 0$.
We have a quasi-isomorphism
$N^*(Y)\mathop{\oslash}^L N_*(X)\,\rightarrow\,N^*(\Map(X,Y))$.}

\smallskip
As in the classical situation of Lannes' $T$ functor,
we can introduce profinite structures
in order to remove the finiteness assumption ({\it cf}. [\bibref{Mrl}]).

\tit{\it The comparison morphism}\label{ComparisonMorphism}\endtit
We define a natural morphism $N^*(Y)\mathop{\oslash}^L N_*(X)\,\rightarrow\,N^*(\Map(X,Y))$.

An $\E$-algebra $A$ has a universal cofibrant resolution $F_A\,\rto{\sim}{}\,A$ ({\it cf}. [\bibref{GeJ}]).
We recall the definition of this cofibrant algebra in section \ref{AlmostFree}.
We consider the universal cofibrant resolution $F_X = F_{N^*(X)}$
of the cochain algebra of a simplicial set $A = N^*(X)$.

By functoriality,
the evaluation map $X\times\Map(X,Y)\,\rightarrow\,Y$
gives rise to a morphism of $\E$-algebras $F_Y\,\rightarrow\,F_{X\times\Map(X,Y)}$.
Hence,
we have a commutative diagram
$$\xymatrix{ F_Y\ar[d]_{\sim}\ar[r] & F_{X\times\Map(X,Y)}\ar[d]_{\sim} & \\
N^*(Y)\ar[r] & N^*(X\times\Map(X,Y)) \\ }$$
The Eilenberg-Zilber equivalence provides a morphism of dg-modules (the shuffle morphism):
$$N^*(X\times\Map(X,Y))\,\rightarrow\,N^*(X)\otimes N^*(\Map(X,Y)).$$
But, the shuffle morphism is not a morphism of $\E$-algebras.
Nevertheless,
according to the next theorem,
we have a morphism of $\E$-algebras
$$F_{X\times\Map(X,Y)}\,\rightarrow\,N^*(X)\otimes N^*(\Map(X,Y))$$
which is cohomologically equivalent to the shuffle morphism
(we refer to theorem \ref{EZKunnethEquivalence} for a more precise statement).

We assume that $X$ is a finite simplicial set so that $\Hom_\F(N_*(X),-) = N^*(X)\otimes -$.
We consider the adjoint of the composite morphism
$$\displaylines{ F_Y\,\rightarrow\,F_{X\times\Map(X,Y)}\,\rightarrow\,N^*(X)\otimes N^*(\Map(X,Y)). \cr
\noalign{\hbox{This morphism}}
F_Y\oslash N_*(X)\,\rightarrow\,N^*(\Map(X,Y)) \cr }$$
is functorial in $X$ and $Y$.
Hence,
we are done in this case.

We observe that both functors
$X\mapsto F_Y\oslash N_*(X)$
and
$X\mapsto N^*(\Map(X,Y))$
preserves colimits.
Therefore,
the comparison morphism is also defined for $X$ infinite
by a colimit argument.
Furthermore,
the general case of theorem \ref{resolutionMapSpace}
is implied by the case of a finite simplicial set $X$.

\th\label{EZequivalence}\endtit
{\it The shuffle morphism $N^*(X\times Y)\,\rightarrow\,N^*(X)\otimes N^*(Y)$
is not a morphism of $\E$-algebras,
but it extends to a morphism of $\E$-algebras
$F_{X\times Y}\,\rightarrow\,N^*(X)\otimes N^*(Y)$
functorially in $X$ and $Y$.}

\tit{\it On pointed simplicial sets and augmented algebras}\label{ReducedAlgebras}\endtit
We have similar statements for pointed spaces.
In this context,
we consider the reduced Barratt-Eccles operad $\tilde{\E}$
which is defined by the relation
$$\tilde{\E}(r) = \left\{\matrix{ 0,\hfill & \ \hbox{if}\ r=0,\hfill\cr
\E(r),\hfill &\ \hbox{otherwise}.\hfill\cr }\right.$$
An algebra over the reduced Barratt-Eccles operad
is equivalent to the augmentation ideal
of an augmented $\E$-algebra.
Precisely,
an augmented $\E$-algebra $A$
is equipped with a morphism of $\E$-algebras $\epsilon: A\,\rightarrow\,\F$.
Since the augmentation $\epsilon: A\,\rightarrow\,\F$ is supposed to preserve the unit morphisms,
we have a splitting $A = \F\oplus\tilde{A}$,
where $\tilde{A} = \ker(\epsilon: A\,\rightarrow\,\F)$
is the augmentation ideal of $A$.
The augmentation ideal $\tilde{A}$ is preserved by the operations $\rho: A^{\otimes r}\,\rightarrow\,A$ where $r\geq 1$.
Therefore, the module $\tilde{A}$ is equipped with the structure of an $\tilde{\E}$-algebra.
Conversely,
the evaluation products $\E(r)\otimes A^{\otimes r}\,\rightarrow\,A$
are determined by the algebra structure of the augmentation ideal $\tilde{A}$,
because the evaluation of an operation at an algebra unit
$\rho(a_1,\ldots,\eta(1),\ldots,a_r)\in A$
agrees with the action of the composite operation
$\rho(1,\ldots,\eta,\ldots,1)\in\E(r-1)$
(see [\bibref{F}, Section 1] for more details).
We have dual observations for coalgebras.

For a simplicial set $X$,
the choice of a basepoint $\pt\,\rightarrow\,X$
determines a morphism of $\E$-algebras
$N^*(X)\,\rightarrow\,N^*(\pt) = \F$.
Hence,
the cochain complex of a pointed simplicial set is an augmented $\E$-algebra.
The reduced cochain complex $\tilde{N}^*(X)$ is the augmentation ideal of $N^*(X)$.

\th\label{resolutionPointedMapSpace}\endtit
{\it Let $X$ and $Y$ be pointed simplicial sets.
We assume that $\pi_n(Y)$ is a finite $p$-group for $n\geq 0$.
We have a quasi-isomorphism
$\tilde{N}^*(Y)\mathop{\oslash}^L\tilde{N}_*(X)\,\rightarrow\,\tilde{N}^*(\Map_*(X,Y))$.}

\tit{\it The division functor in the context of the reduced Barratt-Eccles operad}\label{ReducedDivision}\endtit
We assume that $\tilde{F}$ is the augmentation ideal of a cofibrant $\E$-algebra $F = \F\oplus\tilde{F}$.
We assume that $\tilde{K}$ is the coaugmentation coideal of an $\E$-coalgebra $K = \F\oplus\tilde{K}$.
In this context,
the construction of paragraph \ref{TensorAlgebra}
provides a functor in the category of algebras over the reduced Barratt-Eccles operad
$\Hom_\F(\tilde{K},-): \tilde{\E}\Alg\,\rightarrow\,\tilde{\E}\Alg$.
In the theorem,
we consider the associated adjoint functor
$-\oslash\tilde{K}: \tilde{\E}\Alg\,\rightarrow\,\tilde{\E}\Alg$.

Let us introduce the augmented algebra $F\oslash\tilde{K}\in\E\Alg$
associated to $\tilde{F}\oslash\tilde{K}\in\tilde{\E}\Alg$
and such that $F\oslash\tilde{K} = \F\oplus\tilde{F}\oslash\tilde{K}$.
We have a cocartesian square of augmented algebras
$$\xymatrix{ F\ar[r]\ar[d] & F\oslash K\ar[d] \\
\F\ar[r] & F\oslash\tilde{K} \\ }$$
as proved by an easy inspection.
If we assume that $F = F_Y$ is a cofibrant model of $N^*(Y)$,
then we conclude from this observation
that theorem \ref{resolutionPointedMapSpace} is a consequence of theorem \ref{resolutionMapSpace}.
More precisely,
the morphism
$F_Y\oslash\tilde{N}_*(X)\,\rightarrow\,N^*(\Map_*(X,Y))$
is a quasi-isomorphism,
because the cofibration sequence $F_Y\,\rightarrow\,F_Y\oslash N_*(X)\,\rightarrow\,F_Y\oslash\tilde{N}_*(X)$
is a model of the fibration sequence $\Map_*(X,Y)\,\rightarrow\,\Map(X,Y)\,\rightarrow\,Y$
({\it cf}. [\bibref{Mnd}, Section 5]).

\tit{\it Mandell's comparison result}\endtit
The map $X\,\mapsto\,N^*(X)$ determines a contravariant functor $N^*(-): \Spaces^\op\,\rightarrow\,\E\Alg$.
This functor has a left adjoint $G: \E\Alg\,\rightarrow\,\Spaces^\op$.
These functors determine a pair of derived adjoint functors
$$L G: \Ho(\E\Alg)\,\rightarrow\,\Ho(\Spaces)^\op
\qquad\hbox{and}
\qquad N^*(-): \Ho(\Spaces)^\op\,\rightarrow\,\Ho(\E\Alg).$$

To be explicit,
given any $\E$-algebra $A$,
we set $G(A) = \Hom_\E(A,N^*(\Delta^\bullet))$.
The adjunction relation reads
$$\Hom_\Spaces(X,G(A)) = \Hom_\E(A,N^*(X)).$$
We have also $L G(A) = \Hom_\E(F,N^*(\Delta^\bullet))$,
where $F$ is a cofibrant resolution of $A$.

A simplicial set $X$ is resolvable by $\E$-algebras
if the adjonction unit
$X\,\rightarrow\,L G(N^*(X))$
is a weak-equivalence.
According to Mandell ({\it cf}. [\bibref{Mnd}]),
this property holds if $X$ is connected, nilpotent, $p$-complete
and has a finite $p$-type
(provided the ground field is the algebraic closure $\F = \bar{\F}_p$).

\smallskip
In our framework, we have the following result:

\th\endtit{\it The mapping spaces $\Map(X,Y)$,
where $X$ and $Y$ are as in theorem \ref{resolutionPointedMapSpace},
are resolvable by $\E$-algebras.}

\smallskip
We have a canonical map $\Map(X,Y)\,\rightarrow\,\Hom_\E(F_Y,N^*(X)\otimes N^*(\Delta^\bullet))$.
This map is equivalent to a morphism $\Map(X,Y)\,\rightarrow\,\Hom_\E(F_Y\oslash N_*(X),N^*(\Delta^\bullet))$.
To be precise,
an $n$-dimensional morphism $u\in\Map(X,Y)_n$
is a map $u: X\times\Delta^n\rightarrow Y$
and determines a morphism of $\E$-algebras
$$F_Y\,\rightarrow\,F_{X\times\Delta^n}$$
which can be composed with the morphism
$$F_{X\times\Delta^n}\,\rightarrow\,N^*(X)\otimes N^*(\Delta^n)$$
of theorem \ref{EZequivalence}.

It is straighforward to extend the arguments of Mandell
and to prove that this canonical map
$\Map(X,Y)\,\rightarrow\,\Hom_\E(F_Y,N^*(X)\otimes N^*(\Delta^\bullet))$
is a weak equivalence.
This assertion is also related to an unpublished work of Dwyer and Hopkins.

\section{Cofibrant resolutions and division functors}

\subsection{Free algebras and division functors}

\tit{\it Free algebras}\endtit
The free $\E$-algebra generated by an $\F$-module $V$ is denoted by $\E(V)$.
Let us remind the reader that a realization of $\E(V)$
is provided by a generalization of the symmetric algebra
({\it cf}. [\bibref{GeJ}], [\bibref{GiK}]).
To be explicit, we have
$$\E(V) = \bigoplus_{r=0}^\infty\E_{(r)}(V)
\qquad\hbox{where}\qquad\E_{(r)}(V) = (\E(r)\otimes V^{\otimes r})_{\Sigma_r}.$$
Hence, a tensor $\rho\otimes(v_1\otimes\cdots\otimes v_r)\in\E(r)\otimes V^{\otimes r}$
gives rise to an element of $\E(V)$
which we denote by $\rho(v_1,\ldots,v_r)\in\E(V)$.
By construction,
we have the relation $(w\rho)(v_1,\ldots,v_r) = \pm\rho(v_{w(1)},\ldots,v_{w(r)})$.
Observe that $V = \E_{(1)}(V)\subset\E(V)$,
because an element $v\in V$ can be identified with the tensor $1(v)\in\E(V)$,
where $1\in\E(1)$ is the unit of the Barratt-Eccles operad.

The evaluation product of the free $\E$-algebra
$\E(r)\otimes\E(V)^{\otimes r}\,\rightarrow\,\E(V)$
is induced by the composition product of the Barratt-Eccles operad
$\E(r)\otimes\E(s_1)\otimes\cdots\otimes\E(s_r)\,\rightarrow\,\E(s_1+\cdots+s_r)$.
One should observe that the element $\rho(v_1,\ldots,v_r)\in\E(V)$
represents also the image of $v_1,\ldots,v_r\in V$
under the operation $\rho: \E(V)^{\otimes r}\,\rightarrow\,\E(V)$
(this follows from the unit relation in the Barratt-Eccles operad).

The free $\E$-algebra is characterized by the adjunction relation $\Hom_\E(\E(V),A) = \Hom_\F(V,A)$.
We make the $\E$-algebra morphism $\phi_f: \E(V)\,\rightarrow\,A$
associated to an $\F$-module morphism $f: V\,\rightarrow\,A$
explicit in the paragraph \ref{FreeMorphism} below.
Let us recall that a colimit of free algebras is a free algebra.
More precisely,
because of the adjunction relation,
we obtain $\colim_i \E(V_i) = \E(\colim_i V_i)$.

\smallskip
We observe that the adjunction relations
$\Hom_\E(\E(V\otimes K),B) = \Hom_\F(V\otimes K,B) = \Hom_\F(V,\Hom_\F(K,B)) = Hom_\E(\E(V),\Hom_\F(K,B))$
gives immediately:

\lem\endtit{\it We have $\E(V)\oslash K = \E(V\otimes K)$.}

\tit{\it Free algebra morphisms}\label{FreeMorphism}\endtit
By the universal property of a free algebra,
we assume that a morphism of dg-modules $f: V\,\rightarrow\,A$,
where $A$ is an $\E$-algebra,
extends to one and only one morphism of $\E$-algebras,
which we denote by $\phi_f: \E(V)\,\rightarrow\,A$.
In fact,
if $\phi: \E(V)\,\rightarrow\,A$ is a morphism of $\E$-algebras,
then we have the relation
$\phi(\rho(v_1,\ldots,v_r)) = \rho(\phi(v_1),\ldots,\phi(v_r))$.
Therefore,
a morphism of $\E$-algebras $\phi: \E(V)\,\rightarrow\,A$ is determined by its restriction to $V\subset\E(V)$
and is equivalent to a morphism of $\F$-modules $f: V\,\rightarrow\,A$.
Furthermore,
the algebra morphism associated to $f: V\,\rightarrow\,A$
can be defined explicitly by the formula
$$\phi_f(\rho(v_1,\ldots,v_r)) = \rho(f(v_1),\ldots,f(v_r)).$$
(Hence, for $v\in V$, we have $\phi_f(v) = f(v)$.)

\smallskip
We make explicit the morphism $\phi_f\oslash K: \E(V\otimes K)\,\rightarrow\,\E(W\otimes K)$
associated to a morphism of free algebras $\phi_f: \E(V)\,\rightarrow\,\E(W)$.
We denote also the morphism of dg-modules
which is equivalent to $\phi_f\oslash K: \E(V\otimes K)\,\rightarrow\,\E(W\otimes K)$
by the notation $f\oslash K: V\otimes K\,\rightarrow\,\E(W\otimes K)$.

\lem\label{DivisionFreeMorphism}\endtit{\it The image of an element $v\otimes c\in V\otimes K$
under the morphism $\phi_f\oslash K: \E(V\otimes K)\,\rightarrow\,\E(W\otimes K)$
is determined as follows. We assume
$$f(v) = \sum\nolimits_i \rho^i(w_{(1)}^i,\ldots,w_{(r)}^i)\in\E(W).$$
We consider the diagonal $\Delta(\rho^i) = \sum_j \rho^j_{(1)}\otimes\rho^j_{(2)}$
of the operations $\rho^i\in\E(r)$ which occur in the expansion above.
We take the image of $c\in K$ under the cooperations $\rho^j_{(1)}{}^*: K\,\rightarrow\,K^{\otimes r}$:
$$\rho^j_{(1)}{}^*(c)
= \sum\nolimits_k c_{(1)}^k\otimes\cdots\otimes c_{(r)}^k\in K^{(r)}.$$
We have finally:
$$(f\oslash K)(v\otimes c)
= \sum\nolimits_{i j k} \pm\rho^j_{(2)}(w_{(1)}^i\otimes c_{(1)}^k,\ldots,w_{(r)}^i\otimes c_{(r)}^k)
\in\E(W\otimes K).$$}

\smallskip
The verification of this assertion is straightforward
from the explicit form of the adjunction relation $\Hom_\E(\E(V\otimes K),A) = \Hom_\E(\E(V),\Hom_\F(K,A))$
and from the definition of the $\E$-algebra structure of $\Hom_\F(K,A)$.

\subsection{Cell algebras and division functors}

\tit{\it Cell extensions}\endtit
We consider the suspension sequence $\Sigma^* V\,\rightarrow\,C^* V\,\rightarrow\,V$
of an upper graded differential module $V$.
We have $C^* V = \F\,e\otimes V\oplus\F\,b\otimes V$,
where $\deg(e) = 1$, $\deg(b) = 0$ and $\delta(b) = e$.
The suspension $\Sigma^* V$ is identified with the module $\F\,e\otimes V\subset C^* V$.
Similarly, we have an isomorphism $V\simeq\F\,b\otimes V$
and the module $V$ can be identified with a quotient
of $C^* V$.
We have a morphism of $\E$-algebras $\E(\Sigma^* V)\,\rightarrow\,\E(C^* V)$
induced by the canonical morphism $\Sigma^* V\,\rightarrow\,C^* V$.

In the context of algebras over the Barratt-Eccles operad,
a cell extension is an $\E$-algebra morphism $A\,\rightarrow\,A\vee_f\E(C^*V)$
provided by a cocartesian square
$$\xymatrix{ \E(\Sigma^* V)\ar[r]\ar[d] & \E(C^* V)\ar[d] \\ A\ar[r] & A\vee_f\E(C^*V) \\ }$$
where the left-hand side arrow is determined by a morphism of dg-modules $f: \Sigma^* V\,\rightarrow\,A$
(the analogue of an attaching map).

\smallskip
The next assertion is a general property of left adjoint functors:

\lem\label{DivisionCellExtension}\endtit{\it The division functor preserves cell extensions.
We have explicitly $\bigl(A\vee_f \E(C^* V)\bigr)\oslash K = (A\oslash K)\vee_{f\oslash K}\E(C^* V\otimes K)$.}

\tit{\it Cell algebras}\endtit
A cell algebra $F$ is a colimit of a sequence of cell extensions
$$\F = F^{(-1)}\,\rightarrow\,F^{(0)}\,\rightarrow\,\cdots\,\rightarrow\,F^{(n)}\,\rightarrow\,\cdots$$
Hence,
we have explicitly $F^{(n)} = F^{(n-1)}\vee_{f^{(n)}}\E(C^* V^{(n)})$
where $f^{(n)}: \Sigma^* V^{(n)}\,\rightarrow\,F^{(n-1)}$.

A cell algebra is a cofibrant object in the category of $\E$-algebras.
In fact, it is known that the cofibrant objects in the category of $\E$-algebras
can be characterized as retracts of cell algebras ({\it cf}. [\bibref{Mnd}, Section 2]).
Furthermore,
any $\E$-algebra $A$ has a cofibrant resolution $F\,\rto{\sim}{}\,A$
where $F$ is a cell algebra.

More generally,
a morphism $A\,\rightarrow\,F$
is a relative cell inclusion
if we have a sequence of cell extensions
$$A = F^{(-1)}\,\rightarrow\,F^{(0)}\,\rightarrow\,\cdots\,\rightarrow\,F^{(n)}\,\rightarrow\,\cdots$$
such that $F = \colim_n F^{(n)}$.

\smallskip
The next assertion is an immediate consequence of lemma \ref{DivisionCellExtension}:

\lem\label{DivisionCellAlgebra}\endtit{\it The division functor preserves cell structures.
More explicitly,
if $F$ is a cell algebra, as in the paragraph above,
then we have $F\oslash K = \colim_n F^{(n)}\oslash K$,
where $F^{(n)}\oslash K = F^{(n-1)}\oslash K\vee_{f^{(n)}\oslash K}\E(C^* V^{(n)}\otimes K)$.}

\smallskip
This lemma allows to determine the image of an $\E$-algebra $A\in\E\Alg$
under the derived division functor $-\mathop{\oslash}^L K: \Ho(\E\Alg)\,\rightarrow\,\Ho(\E\Alg)$.
Just recall that a representative of $A\mathop{\oslash}^L K\in\E\Alg$
is given by an $\E$-algebra $F\oslash K\in\E\Alg$,
where $F$ is a cofibrant resolution of $A$.

\subsection{Almost free algebras and division functors}\label{AlmostFree}

\tit{\it Free algebra derivations}\endtit
A homogeneous morphism $d: A\,\rightarrow\,A$,
where $A$ is an $\E$-algebra,
is a derivation
if we have the relation
$d(\rho(a_1,\ldots,a_r)) = \sum_{i=1}^r\pm\rho(a_1,\ldots,d(a_i),\ldots,a_r)$.
In the case of a free $\E$-algebra $A = \E(V)$,
a homogeneous morphism $h: V\,\rightarrow\,\E(V)$
extends to one and only one derivation $d_h: \E(V)\,\rightarrow\,\E(V)$.
In fact,
as for free algebra morphisms,
a derivation $d: \E(V)\,\rightarrow\,\E(V)$ is determined by its restriction to $V\subset\E(V)$,
because of the derivation relation.
Furthermore,
the derivation $d_h: \E(V)\,\rightarrow\,\E(V)$
can be defined explicitly by the formula
$$d_h(\rho(v_1,\ldots,v_r)) = \sum_{i=1}^r\pm\rho(v_1,\ldots,h(v_i),\ldots,v_r).$$
(Hence, for $v\in V$, we have $d_h(v) = h(v)$.)

\tit{\it Almost free algebras}\label{AlmostFreeAlgebra}\endtit
The free $\E$-algebra $F = \E(V)$
is equipped with the differential $\delta: \E(V)\,\rightarrow\,\E(V)$
induced by the differential of $V$
and by the differential of the Barratt-Eccles operad $\E$.
Explicitly,
the differential of an element $\rho(v_1,\ldots,v_r)\in\E(V)$
is given by the formula
$$\delta(\rho(v_1,\ldots,v_r)) = (\delta\rho)(v_1,\ldots,v_r)
+ \sum_{i=1}^r\pm\rho(v_1,\ldots,\delta(v_i),\ldots,v_r).$$

An almost free $\E$-algebra $F$
is a free $\E$-algebra $F = \E(V)$
equipped with a differential
which differs from the canonical one by a derivation.
To be precise,
an almost free $\E$-algebra $F = \E(V)$
is equipped with a differential $\delta_h: \E(V)\,\rightarrow\,\E(V)$ such that $\delta_h = \delta + d_h$,
where $d_h: \E(V)\,\rightarrow\,\E(V)$ is a derivation associated to a homogeneous morphism of degree $1$.
The differential $\delta_0: \E(V)\,\rightarrow\,\E(V)$
associated to the $0$ morphism
reduces to the canonical differential.

\tit{\it Morphisms of almost free algebras}\endtit
We consider a morphism $\phi: F\,\rightarrow\,A$ where $F$ is an almost free algebra $F = \E(V)$.
In this situation,
we assume the relation $\delta\phi = \phi\delta_h$,
where $\delta_h: \E(V)\,\rightarrow\,\E(V)$ denotes the differential of $F$.
We have necessarily $\phi = \phi_f$
where $f: V\,\rightarrow\,A$ is a homogeneous morphism of degree $0$.
But, in the case of an almost free algebra,
the morphism $f: V\,\rightarrow\,A$ is not supposed to preserve the differentials of dg-modules,
because the relation $\delta\phi_f = \phi\delta_h$
does not depend on this property.

\tit{\it The strict category of almost free algebras}\label{StrictCategory}\endtit
We consider the category whose objects are almost free algebras
and whose morphisms are strict morphisms of almost free algebras.
By definition, a morphism of almost free algebras $\phi_f: \E(V)\,\rightarrow\,\E(W)$ is strict
if we have $\phi_f(V)\subset W$
(equivalently,
if the morphism $\phi_f: \E(V)\,\rightarrow\,\E(W)$
is induced by a morphism of dg-modules $f: V\,\rightarrow\,W$).

We have an obvious forgetful functor from the strict category of almost free $\E$-algebras
to the category of all $\E$-algebras.
We claim that this forgetful functor has a right adjoint.
To be more explicit,
to any $\E$-algebra $A\in\E\Alg$,
we associate an almost free algebra $F_A = \E(C_A)$
together with a morphism $F_A\,\rightarrow\,A$.
Furthermore,
any morphism $F\,\rightarrow\,A$,
where $F = \E(V)$ is almost free,
has a unique factorization
$$\xymatrix{ & F_A\ar[d] \\ F\ar@{-->}[ur]^{\tilde{\phi}}\ar[r] & A \\ }$$
such that $\tilde{\phi}: F\,\rightarrow\,F_A$ is a strict morphism of almost free algebras.

In fact, the augmentation morphism $F_A\,\rightarrow\,A$ is a quasi-isomorphism,
and we explain in section \ref{CellFreeStructure}
that the $\E$-algebra $F_A$
is a canonical cell resolution of $A$
({\it cf}. proposition \ref{UniversalCellResolution}).
The purpose of the next paragraphs is to recall the construction of the almost free algebra $F_A$
in connection with the operadic bar construction of Getzler-Jones ({\it cf}. [\bibref{GeJ}])
and Ginzburg-Kapranov ({\it cf}. [\bibref{GiK}]).

\tit{\it The structure of the bar cooperad}\endtit
To be precise,
we consider the bar cooperad $\B\E$ introduced by E. Getzler and J. Jones ({\it cf}. [\bibref{GeJ}]).
We recall that the structure of a cooperad consists of a sequence $\B\E(r)$, $r\in\N$,
where $\B\E(r)$ is a representation of the symmetric group $\Sigma_r$;
together with decomposition coproducts
$$\B\E(n)\,\rightarrow\,\B\E(r)\otimes\B\E(s_1)\otimes\cdots\otimes\B\E(s_r),$$
defined for $n\in\N$ and $s_1+\cdots+s_r=n$;
but, the structure is such that for a fixed element of $\B\E(n)$
only finitely many of these coproducts are non zero.
We assume the dual equivariance and associativity properties of an operad composition product.
Therefore,
we may also consider the operad $\B\E^{\vee}$ formed by the dual representations $\B\E(r)^{\vee}$, $r\in\N$,
of the modules $\B\E(r)$, $r\in\N$.
(But, these structures are not equivalent.)

Similarly,
the structure of a $\B\E$-coalgebra consists of a dg-module $V$
equipped with a coproduct $V\,\rightarrow\,\bigoplus_{r=0}^\infty(\B\E(r)\otimes V^{\otimes r})_{\Sigma_r}$.
(We observe that the module of coinvariants $(\B\E(r)\otimes V^{\otimes r})_{\Sigma_r}$
is isomorphic to a module of invariants $(\B\E(r)\otimes V^{\otimes r})^{\Sigma_r}$,
because,
in the case of the Barratt-Eccles operad,
the module of coefficients $\B\E(r)$ is a regular representation of $\Sigma_r$.)
Equivalently,
we have a coproduct
$V\,\rightarrow\,\B\E(r)\otimes V^{\otimes r}$
for $r\in\N$;
but, the structure is such that for a fixed element of $V$
only finitely many of these coproducts are non zero.
Finally,
a coalgebra over a cooperad $\B\E$
is equivalent a coalgebra over a dual operad $\B\E^{\vee}$
together with a certain nilpotence assumption.
The construction
$$\B\E(V) = \bigoplus_{r=0}^\infty\B\E_{(r)}(V),
\qquad\hbox{where}
\qquad\B\E_{(r)}(V) = (\B\E(r)\otimes V^{\otimes r})_{\Sigma_r},$$
provides a realization of the cofree $\B\E$-coalgebra.

\tit{\it The structure of almost free algebras}\endtit
The bar cooperad $\B\E$ is characterized by the following properties.
The structure of a $\B\E$-coalgebra
$\B\E(r)\otimes V\,\rightarrow\,V^{\otimes r}$
is equivalent to a differential
$\delta_h: \E(V)\,\rightarrow\,\E(V)$.
Furthermore,
a morphism of $\B\E$-coalgebras
is equivalent to a strict morphism of almost free algebras
$\phi_f: \E(V)\,\rightarrow\,\E(W)$.
In fact,
a morphism $f: V\,\rightarrow\,W$ preserves $\B\E$-coalgebra structures
if and only if the associated morphism of free algebras
$\phi_f: \E(V)\,\rightarrow\,\E(W)$
preserves the differentials of almost free algebras.

The almost free algebra $F_A = \E(C_A)$ is associated to a $\B\E$-coalgebra $C_A$ which has the following structure.
The cofree $\B\E$-coalgebra $\B\E(A)$
is equipped with a differential $\delta: \B\E(A)\,\rightarrow\,\B\E(A)$
induced by the internal differential of $A$ and by the differential of $\B\E$.
As in the context of almost free algebras over an operad,
we consider a differential $\delta_A: \B\E(A)\,\rightarrow\,\B\E(A)$
which differs from $\delta: \B\E(A)\,\rightarrow\,\B\E(A)$
by a specific coderivation of the cofree cooperad $\B\E(A)$.
(In fact, the coderivation is determined by the algebra structure of $A$.)
The coalgebra $C_A$ consists of the cofree coalgebra $\B\E(A)$
equipped with this differential $\delta_A: \B\E(A)\,\rightarrow\,\B\E(A)$.

\tit{\it Almost free algebras and division functors}\endtit
We can make explicit the image of an almost free algebra under a division functor (as for cell algebras).
To be precise,
we observe that the division functor preserves almost free algebras.
In fact,
if $F = \E(V)$,
then we have $F\oslash K = \E(V\otimes K)$.
It is also straighforward to determine the derivation
$d_h\oslash K: \E(V\otimes K)\,\rightarrow\,\E(V\otimes K)$
which yields the differential of $F\oslash K$.
The morphism of dg-modules
which determines this derivation
is also denoted by $h\oslash K: V\otimes K\,\rightarrow\,\E(V\otimes K)$.
In fact, we obtain the same result as in lemma \ref{DivisionFreeMorphism}:

\lem\label{AlmostFreeDivision}\endtit{\it The image of an element $v\otimes c\in V\otimes K$
under the derivation $d_h\oslash K: \E(V\otimes K)\,\rightarrow\,\E(V\otimes K)$
is determined as follows. We assume
$$h(v) = \sum\nolimits_i \rho^i(v_{(1)}^i,\ldots,v_{(r)}^i)\in\E(V).$$
We consider the diagonal $\Delta(\rho^i) = \sum_j \rho^j_{(1)}\otimes\rho^j_{(2)}$
of the operations $\rho^i\in\E(r)$ which occur in the expansion above.
We take the image of $c\in K$ under the cooperations $\rho^j_{(1)}{}^*: K\,\rightarrow\,K^{\otimes r}$:
$$\rho^j_{(1)}{}^*(c)
= \sum\nolimits_k c_{(1)}^k\otimes\cdots\otimes c_{(r)}^k\in K^{(r)}.$$
We have finally:
$$(h\oslash K)(v\otimes c)
= \sum\nolimits_{i j k} \pm\rho^j_{(2)}(v_{(1)}^i\otimes c_{(1)}^k,\ldots,v_{(r)}^i\otimes c_{(r)}^k)
\in\E(V\otimes K).$$}

\subsection{Cell algebras are almost free}\label{CellFreeStructure}

\tit{\it Cell extensions of almost free algebras}\label{CellExtensionAlmostFree}\endtit
Let us consider a cell extension $F\,\rightarrow\,F\vee_f\E(C^* V)$,
where $F = \E(U)$ is almost free
and is equipped with the differential $\delta_h: \E(U)\,\rightarrow\,\E(U)$.
In this case,
the algebra $F\vee_f\E(C^* V)$ is also almost free.
Precisely,
we obtain $F\vee_f\E(C^* V) = \E(U\oplus V)$,
because the free algebra functor preserves colimits.
Furthermore,
the morphism $h': U\oplus V\,\rightarrow\,\E(U\oplus V)$
which determines the differential of $F\vee_f\E(C^* V)$
maps the module $V\subset U\oplus V$ into $\E(U)\subset\E(U\oplus V)$.
In fact,
this map is determined by the derivation of $F$ on $U\subset U\oplus V$
and by the attaching morphism on $V\subset U\oplus V$.
Precisely,
the map $h': U\oplus V\,\rightarrow\,\E(U\oplus V)$
is the sum of the composites
$$U\,\rto{h}{}\,\E(U)\,\hookrightarrow\,\E(U\oplus V)
\qquad\hbox{and}
\qquad V\,\simeq\,\Sigma^* V\,\rto{f}{}\,\E(U)\,\hookrightarrow\,\E(U\oplus V).$$

By induction,
we deduce from these observations that a cell algebra is always an almost free algebra.
Conversely,
we have the following result:

\prop\label{AlmostFreeCellAlgebras}\endtit
{\it A cell algebra is equivalent to an almost free algebra $F' = \E(U')$,
where the dg-module $U'$ has a filtration
$$0 = U^{(-1)}\subset U^{(0)}\subset\cdots\subset U^{(n)}\subset\cdots\subset U'$$
such that the morphism $h': U'\,\rightarrow\,\E(U')$,
which determines the differential of $F' = \E(U')$,
verifies $h'(U^{(n)})\subset\E(U^{(n-1)})$.}

\proof\endtit
Suppose given an almost free algebra $F' = \E(U')$ as in the proposition above.
We consider the algebra $F^{(n)} = \E(U^{(n)})$
together with the differential determined by the restriction of $h': U'\,\rightarrow\,\E(U')$
to $U^{(n)}\subset U'$.
We fix a splitting $U^{(n)} = U^{(n-1)}\oplus V^{(n)}$.
We consider the restriction of $h: U'\,\rightarrow\,\E(U')$
to $V^{(n)}\subset U'$.
According to the observations of paragraph \ref{CellExtensionAlmostFree},
this gives a morphism of dg-modules $f^{(n)}: \Sigma^* V^{(n)}\,\rightarrow\,F^{(n-1)}$
such that $F^{(n)} = F^{(n-1)}\vee_{f^{(n)}}\E(C^* V^{(n)})$.

\smallskip
Proposition \ref{AlmostFreeCellAlgebras} can be generalized to a relative context.
Namely, for an almost free algebra $F = \E(U)$,
a relative cell inclusion $F\,\rightarrow\,F'$
is equivalent to a strict morphism $\phi_f: \E(U)\,\rightarrow\,\E(U')$,
where the dg-module $V$ has a filtration
$$U = U^{(-1)}\subset U^{(0)}\subset\cdots\subset U^{(n)}\subset\cdots\subset U'$$
such that the map $h': U'\,\rightarrow\,\E(U')$,
which determines the differential of $F' = \E(U')$,
verifies $h'(U^{(n)})\subset\E(U^{(n-1)})$.

\lem\label{CellStrictMorphism}\endtit{\it Suppose given a strict morphism $\phi_f: F\,\rightarrow\,F'$
of almost free algebras $F = \E(U)$ and $F' = \E(U')$
induced by an injective morphism of dg-modules $f: U\,\rightarrow\,U'$.
If the almost free algebra $F' = \E(U')$ is equipped with a cell structure,
then this morphism $\phi_f: F\,\rightarrow\,F'$
is equivalent to a relative cell inclusion (and is a cofibration).}

\proof\endtit
As in proposition \ref{AlmostFreeCellAlgebras},
we let $h': U'\,\rightarrow\,\E(U')$
denote the map which determines the differential of $F' = \E(U')$.
The cell structure of $F' = \E(U')$
is determined by a filtration
such that
$$0 = U^{(-1)}\subset U^{(0)}\subset\cdots\subset U^{(n)}\subset\cdots\subset U'.$$
We can assume that $U$ is a sub-dg-module of $U'$
and $F = \E(U)$ is a sub-dg-algebra of $F' = \E(U')$.
We have then $h'(U)\subset\E(U)$
so that the differential of $F'$ verifies $d_{h'}(F)\subset F$.
Clearly, the filtration
$$U = U+U^{(-1)}\subset U+U^{(0)}\subset\cdots\subset U+U^{(n)}\subset\cdots\subset U'$$
gives the inclusion morphism $\E(U)\,\rightarrow\,\E(U')$
the structure of a relative cell inclusion.

\smallskip
We consider the almost free algebra $F_A = \E(\B\E(A))$
introduced in paragraph \ref{StrictCategory}.
We prove the following result:

\prop\label{UniversalCellResolution}\endtit
{\it The universal almost free resolution of an $\E$-algebra $F_A = \E(\B\E(A))$ has a cell structure.
Consequently, the construction $A\mapsto F_A$ provides a cofibrant resolution for all $\E$-algebras $A$.}

\proof\endtit
We consider the module $\B\E_n(A)\subset\B\E(A)$ generated by the degree $n$ components of the bar construction
$\B\E_n(A) = \bigoplus_{r=0}^\infty(\B\E_n(r)\otimes A^{\otimes r})_{\Sigma_r}$.
We have an exhaustive filtration
$0 = U^{(-1)}\subset U^{(0)}\subset\cdots\subset U^{(n)}\subset\cdots\subset\B\E(A)$
where $U^{(n)} = \bigoplus_{m\leq n} \B\E_m(A)$.
The internal differential of $A$ maps $\B\E_n(A)$ to $\B\E_n(A)$.
The differential of the bar construction maps $\B\E_n(A)$ to $\E(\bigoplus_{m\leq n-1} \B\E_m(A))$
({\it cf}. [\bibref{GeJ}]).
The proposition follows.

\smallskip
The next statement is an immediate corollary of lemma \ref{CellStrictMorphism}
and proposition \ref{UniversalCellResolution}.

\prop\label{UniversalCellMorphism}\endtit
{\it Let $A\,\rightarrow\,B$ be an injective morphism of $\E$-algebras.
The induced morphism $F_A\,\rightarrow\,F_B$ has a relative cell structure.
Consequently, this morphism is a cofibration in the category of $\E$-algebras.}

\section{The example of loop spaces}

\tit{\it The adjoint functors associated to the circle coalgebra}\endtit
As in paragraph \ref{ReducedAlgebras},
we work in the category of algebras over the reduced Barratt-Eccles operad $\tilde{\E}$
(equivalent to the category of augmented $\E$-algebras).
We replace an augmented $\E$-algebra $A$ by its augmentation ideal $\tilde{A}$,
which is equipped with the structure of an $\tilde{\E}$-algebra.

We consider the reduced chain complex of the circle $\tilde{K} = \tilde{N}_*(S^1)$.
In this case,
the module $\Hom_\F(\tilde{N}_*(S^1),\tilde{B})$
is denoted by $\Sigma^*\tilde{B} = \Hom_\F(\tilde{N}_*(S^1),\tilde{B})$.
and is called the suspension algebra of $\tilde{B}$.
Similarly,
the division algebra $\tilde{A}\oslash\tilde{N}_*(S^1)$
is denoted by $\Omega^*\tilde{A} = \tilde{A}\oslash\tilde{N}_*(S^1)$
amd is called the loop algebra of $A$.
We have the adjunction relation
$$\Hom_{\tilde{\E}}(\Omega^*\tilde{A},\tilde{B}) = \Hom_{\tilde{\E}}(\tilde{A},\Sigma^*\tilde{B}).$$

The next propositions give the homotopy significance of suspension and loop algebra functors.

\prop\endtit{\it We assume $B$ is an $\E$-algebra equipped with an augmentation $B\,\rightarrow\,\F$.
We consider the associated augmentation ideal $\tilde{B}$,
which is an algebra over the reduced Barratt-Eccles operad $\tilde{\E}$.
The suspension algebra $\Sigma^*\tilde{B}$
is a representative of the loop object of $\tilde{B}$
in the homotopy category of $\tilde{\E}$-algebras.
(Hence, this is a suspension object in the opposite category.)}

\proof\endtit
We consider the cone algebra $C^*\tilde{B} = \Hom_\F(\tilde{N}_*(\Delta^1),\tilde{B})$,
which is an acyclic algebra
together with a natural fibration
$C^*\tilde{B}\,\rightarrow\,\tilde{B}$.
The suspension algebra is identified with the kernel of this fibration.
Hence,
we have a cartesian square
$$\xymatrix{ \Sigma^*\tilde{B}\ar[r]\ar[d] & C^*\tilde{B}\ar[d] \\
0\ar[r] & \tilde{B} \\ }$$
and, as a conclusion,
the suspension algebra
satisfies the classical definition of a loop object
in a closed model category ({\it cf}. [\bibref{Q}]).

\prop\label{LoopAlgebra}\endtit
{\it We assume $F$ is a cofibrant $\E$-algebra equipped with an augmentation $F\,\rightarrow\,\F$.
We consider the associated augmentation ideal $\tilde{F}$,
which is a cofibrant algebra over the reduced Barratt-Eccles operad $\tilde{\E}$.
The loop algebra $\Omega^*\tilde{F}$
is a representative of the suspension of $\tilde{F}$
in the homotopy category of $\tilde{\E}$-algebras.
(Hence, this is a loop object in the opposite category.)}

\proof\endtit
We consider the algebra $P^*\tilde{F} = \tilde{F}\oslash\tilde{N}(\Delta^1)$.
According to lemma \ref{PathAlgebra} below,
this algebra is acyclic
and, moreover, we have a natural cofibration $\tilde{F}\,\rightarrow\,P^*\tilde{F}$.
Therefore,
the algebra $P^*\tilde{F}$ deserves to be called a path algebra of $\tilde{F}$.
We observe that the loop algebra $\Omega^*\tilde{F}$ fits in the cocartesian square
$$\xymatrix{ \tilde{F}\ar[r]\ar[d] & P^*\tilde{F}\ar[d] \\
0\ar[r] & \Omega^*\tilde{F} \\ }$$
and the conclusion follows.
Namely,
the loop algebra $\Omega^*\tilde{F}$
satisfies the classical definition of a suspension object
in a closed model category ({\it cf}. [\bibref{Q}]).

\lem\label{PathAlgebra}\endtit{\it As in proposition \ref{LoopAlgebra},
we assume that $\tilde{F}$ is a cofibrant $\tilde{\E}$-algebra.
The canonical morphism $0\,\rightarrow\,\tilde{F}\oslash\tilde{N}_*(\Delta^1)$
is an acyclic cofibration.}

\proof\endtit
The morphism $0\,\rightarrow\,\tilde{F}\oslash\tilde{N}_*(\Delta^1)$ has the left lifting property
with respect to fibrations,
because the diagrams
$$\vcenter{\xymatrix{ & \tilde{A}\ar[d] \\
\tilde{F}\oslash\tilde{N}_*(\Delta^1)\ar@{-->}[ur]\ar[r] & B \\ }}
\qquad\hbox{and}
\qquad\vcenter{\xymatrix{ & \Hom_\F(\tilde{N}_*(\Delta^1),\tilde{A})\ar[d] \\
\tilde{F}\ar@{-->}[ur]\ar[r] & \Hom_\F(\tilde{N}_*(\Delta^1),\tilde{B}) \\ }}$$
are equivalent by adjunction.
To be more precise,
if the morphism $\tilde{A}\,\rightarrow\,\tilde{B}$ is a fibration,
then the induced morphism
$\Hom_\F(\tilde{N}_*(\Delta^1),\tilde{A})\,\rightarrow\,\Hom_\F(\tilde{N}_*(\Delta^1),\tilde{B})$
is still surjective.
Furthermore,
both dg-modules $\Hom_\F(\tilde{N}_*(\Delta^1),\tilde{A})$ and $\Hom_\F(\tilde{N}_*(\Delta^1),\tilde{B})$
are clearly acyclic.
Hence,
the induced morphism
$\Hom_\F(\tilde{N}_*(\Delta^1),\tilde{A})\,\rightarrow\,\Hom_\F(\tilde{N}_*(\Delta^1),\tilde{B})$
is an acyclic fibration.
The lemma follows
since the $\E$-algebra $\tilde{F}$
has the left lifting property
with respect to such morphisms.

\smallskip
We consider the left derived functor $L\Omega^*: \Ho(\tilde{\E}\Alg)\,\rightarrow\,\Ho(\tilde{\E}\Alg)$
such that $L\Omega^* = -\oslash^L\tilde{N}_*(S^1)$.
In the case $X = S^1$,
theorem \ref{resolutionPointedMapSpace} asserts
that the algebra $L\Omega^*\tilde{N}^*(Y) = \tilde{N}^*(Y)\oslash^L\tilde{N}_*(S^1)$
is quasi-isomorphic to the cochain algebra $\tilde{N}^*(\Omega Y) = \tilde{N}^*\Map_*(S^1,Y)$,
provided the homotopy groups $\pi_n(Y)$, $n\in\N$, are finite $p$-groups.
In fact,
in the particular case $X = S^1$,
it is not difficult to state a more precise result:

\th\endtit\label{LoopModel}{\it Let $Y$ be a pointed simplicial set.
We assume that the cohomology modules $H^n(Y)$ are finite dimensional
and that the fundamental group $\pi_1(Y)$ is a finite $p$-group.
We have a quasi-isomorphism $L\Omega^*\tilde{N}^*(Y)\,\rto{\sim}{}\,\tilde{N}^*(\Omega Y)$.}

\proof\endtit
We deduce this result from a theorem of Mandell ({\it cf}. [\bibref{Mnd}, Lemma 5.2]).
We recall the main arguments in order to be careful with the assumptions of the theorem.
We work with augmented $\E$-algebras.

We fix a cofibrant augmented $\E$-algebra $F_Y$
(equivalent to a cofibrant $\tilde{\E}$-algebras $\tilde{F}_Y$)
together with a quasi-isomorphism
$F_Y\,\rto{\sim}{}\,N^*(Y)$.
We consider the cocartesian square of augmented $\E$-algebras
$$\xymatrix{ F_Y\ar[r]\ar[d]\ar@{}[dr]|{(M)} & P^* F_Y\ar[d] \\
\F\ar[r] & \Omega^* F_Y \\ }$$
We compare this diagram to the commutative square of cochain algebras
$$\xymatrix{ N^*(Y)\ar[r]\ar[d]\ar@{}[dr]|{(C)} & N^*(P Y)\ar[d] \\
\F\ar[r] & N^*(\Omega Y) \\ }$$

Precisely,
we have a morphism from the cocartesian square (M) to the commutative square (C). 
In general,
we assume that $F_Y$ is the universal resolution of $N^*(Y)$
and we deduce this assertion from theorem \ref{EZequivalence}.
But,
in the current situation,
we can deduce the existence of fill-in morphisms from straighforward arguments.
Explicitly,
the composite morphism
$F_Y\,\rightarrow\,N^*(Y)\,\rightarrow\,N^*(P Y)$
extends to a morphism of augmented $\E$-algebras $P^* F_Y\,\rightarrow\,N^*(P Y)$,
because the canonical morphism $N^*(P Y)\,\rightarrow\,\F$ is an acyclic fibration,
and the cofibration $F_Y\,\rightarrow\,P^* F_Y$
has the left lifting property
with respect to acyclic fibrations.
We have an induced morphism $\Omega^* F_Y\,\rightarrow\,N^*(\Omega Y)$
because the diagram (M) is cocartesian.
We prove that this morphism is a quasi-isomorphism.

The cocartesian diagram (M)
gives rise to a strongly convergent spectral sequence $E^r_{(M)}\,\Rightarrow\,H^*(\Omega^* F_Y)$,
such that $E^2_{(M)} = \Tor^{H^*(F_Y)}_*(\F,H^*(P^* F_Y)) = \Tor^{H^*(Y)}_*(\F,\F)$
({\it cf}. [\bibref{Mnd}, Corollary 3.6]).
We compare this spectral sequence
with the classical Eilenberg-Moore spectral sequence
$E^r_{(C)}\,\Rightarrow\,H^*(\Omega Y)$.
We have $E^2_{(C)} = \Tor^{H^*(Y)}_*(\F,H^*(P Y)) = \Tor^{H^*(Y)}_*(\F,\F)$.
The assumption of the theorem ensures
the strong convergence of the classical Eilenberg-Moore spectral sequence ({\it cf}. [\bibref{D}]).
The morphism of commutative squares $(M)\,\rightarrow\,(C)$
gives rise to a morphism of spectral sequences $E^r_{(M)}\,\rightarrow\,E^r_{(C)}$.
Therefore,
we obtain an isomorphism $H^*(\Omega^* F_Y)\,\rto{\simeq}{}\,H^*(\Omega Y)$,
because the $E^2$ stage of these spectral sequences agree.
This achieves the proof of theorem \ref{LoopModel}.

\smallskip
The purpose of the next paragraphs is to make the loop algebra functor explicit.
Our result is stated in proposition \ref{LoopConstruction}.
In fact,
we recover a construction of Smirnov ({\it cf}. [\bibref{SmrCLoop}], [\bibref{SmrHLoop}]).

First,
we remind the reader of results about the structure of the chain coalgebra of the circle
({\it cf}. [\bibref{BF}]).

\tit{\it The circle coalgebra}\label{CircleCoalgebra}\endtit
We make explicit the structure of the reduced chain coalgebra $\tilde{K} = \tilde{N}_*(S^1)$
of the circle $S^1 = \Delta^1/\partial\Delta^1$.
We have $\tilde{N}_*(S^1) = \F\,e_1$, where $\deg(e_1) = 1$.

We consider the cocycle $\epsilon_r: \E(r)_*\,\rightarrow\,\F$ introduced in article [\bibref{BF}].
We have explicitly $\epsilon_r(w_0,\ldots,w_d) = \pm 1$,
if the sequence $(w_0(1),\ldots,w_d(1))$ is a permutation of $(1,\ldots,r)$,
and $\epsilon_r(w_0,\ldots,w_d) = 0$, otherwise.
According to this definition,
we assume that the cochain $\epsilon_r: \E(r)_*\,\rightarrow\,\F$ is zero in degree $*\not=r-1$.

For any $\rho\in\E(r)$,
the operation $\rho^*: \tilde{N}_*(S^1)\,\rightarrow\,\tilde{N}_*(S^1)^{\otimes r}$
satisfies the formula $\rho^*(e_1) = \epsilon_r(\rho)\,e_1^{\otimes r}$.

\tit{\it The suspension functor}\label{SuspensionAlgebra}\endtit
We make explicit the structure of the suspension algebra $\Sigma^* A$.
In fact,
we have a suspension functor on dg-modules $V\mapsto \Sigma^* V$
such that $\Sigma^* V$ is deduced from $V$ by a shift of degree.

We set explicitly $\Sigma^* V^* = V^{*-1}$.
We let $\Sigma^* v\in\Sigma^* V$
denote the element associated to $v\in V$.
Clearly,
our definition of the suspension functor in the category of $\tilde{\E}$-algebras
agrees with the definition of the suspension functor for dg-modules.

We consider the cap product with the cocycle $\epsilon_r: \E(r)_*\,\rightarrow\,\F$.
We obtain a morphism of dg-modules $\epsilon_r\cap-: \E(r)_*\,\rightarrow\,\E(r)_{*-r+1}$.
We have explicitly $\epsilon_r\cap\rho = \sum_i \epsilon_r(\rho_{(1)}^i)\,\rho_{(2)}^i$,
where $\Delta(\rho) = \sum_i \rho_{(1)}^i\otimes\rho_{(2)}^i$
is the diagonal of $\rho\in\E(r)$.
We deduce from the construction of paragraph \ref{TensorAlgebra}
that the operation $\rho^*: (\Sigma^* A)^{\otimes r}\,\rightarrow\,\Sigma^* A$
satisfies the relation
$$\rho(\Sigma^* a_1,\ldots,\Sigma^* a_r) = \Sigma^*(\epsilon_r\cap\rho(a_1,\ldots,a_r)).$$

\smallskip
We deduce the structure of a loop algebra $\Omega^* F$ from claim \ref{AlmostFreeDivision}.
We introduce the following definition in order to state our result.

\tit{\it The suspension of the Barratt-Eccles operad}\endtit
As in article [\bibref{BF}],
we consider the suspension of the Barratt-Eccles operad $\Lambda^*\E$
such that $\Lambda^*\E(r)_* = \sgn(r)\otimes\E(r)_{*-r+1}$,
where $\sgn(r)$ denotes the signature representation of the symmetric group $\Sigma_r$.
An algebra over the operad $\Lambda^*\E$
is equivalent to the suspension $\Sigma^* A$ of an $\E$-algebra.
Furthermore,
for free algebras,
we have the relation $\Lambda^*\E(V) = \Sigma^*\E(V^{*+1})$.

The $\E$-algebra structure of a suspension $\Sigma^* A$
can be deduced from the canonical $\Lambda^*\E$-algebra structure
by a restriction process ({\it cf}. [\bibref{BF}]).
Precisely,
the cap product operations $\epsilon_r\cap-: \E(r)_*\,\rightarrow\,\E(r)_{*-r+1}$
define an operad morphism $\epsilon_*\cap-: \E\,\rightarrow\,\Lambda^*\E$.
Consequently,
if $A$ is an $\E$-algebra,
then the canonical evaluation product $\Lambda^*\E(r)\otimes(\Sigma^* A)^{\otimes r}\,\rightarrow\,\Sigma^* A$
restricts to
$$\E(r)\otimes(\Sigma^* A)^{\otimes r}
\,\rightarrow\,\Lambda^*\E(r)\otimes(\Sigma^* A)^{\otimes r}
\,\rightarrow\,\Sigma^* A.$$
We recover clearly the formula of paragraph \ref{SuspensionAlgebra}
from this construction.

\prop\label{LoopConstruction}\endtit{\it If $F$ is an almost free algebra over the Barratt-Eccles operad $\E$,
then $\Sigma^*\Omega^* F$ is an almost free algebra
over the suspension of the Barratt-Eccles operad $\Lambda^*\E$.
More precisely,
if $F = \E(V)$,
then we have $\Sigma^*\Omega^* F = \Lambda^*\E(V)$.
Moreover,
if the differential of $F$ is determined by a morphism
$$V\,\rto{h}{}\,\E(V)$$
(as explained in paragraph \ref{AlmostFreeAlgebra}),
then the almost free algebra $\Sigma^*\Omega^* F$
is equipped with the differential associated to the composite
$$V\,\rto{h}{}\,\E(V)\,\rto{\epsilon_*\cap-}{}\,\Lambda^*\E(V).$$
Hence,
we consider the morphism $\E(V)\,\rto{\epsilon_*\cap-}{}\,\Lambda^*\E(V)$
induced by the cap product with the cocycles $\epsilon_r: \E(r)_*\,\rightarrow\,\F$
of paragraph \ref{CircleCoalgebra}.}

\section{The case of Eilenberg-MacLane spaces}

\smallskip
In this section,
we determine the image of the cofibrant model of an Eilenberg-MacLane space $K(\Z/p,n)$
under a division functor $-\oslash K: \E\Alg\,\rightarrow\,\E\Alg$.
We prove the following result:

\tit{\sc Lemma}\endtit
{\it The canonical morphism $F_{K(\Z/p,n)}\oslash N_*(X)\,\rto{}{}\,N^*\Map(X,K(\Z/p,n))$ is a quasi-isomorphism.}

\smallskip
We replace the universal almost free algebra $F = F_{K(\Z/p,n)}$
by an equivalent cell algebra $F = F_n$
introduced by Mandell.
We prove that we have a quasi-isomorphism
$$F_n\oslash N_*(X)\,\rto{\sim}{}\,N^*\Map(X,K(\Z/p,n))$$
by a generalization of Mandell's arguments.

The general idea consists in introducing Steenrod operations
and unstable algebra structures.
In fact,
the mapping space $\Map(X,K(\Z/p,n))$ is a generalized Eilenberg-MacLane space.
Consequently,
the cohomology algebra $H^*(\Map(X,K(\Z/p,n)),\F)$ is identified with a free object
in the category of unstable algebras over the classical Steenrod algebra.
Our purpose is to identify the cohomology algebra $H^*(F_n\oslash N_*(X))$
with the same free object.

\subsection{The big Steenrod algebra}

\smallskip
We recall the definition and properties of Steenrod operations
in the context of algebras over the Barratt-Eccles operad.
For simplicity,
we assume that the ground field $\F$ has characteristic $p = 2$.
For the general case,
we refer to the article of Mandell ({\it cf}. [\bibref{Mnd}, Sections 11-12]).

\tit{\it Steenrod operations}\endtit
The cohomology $H^*(A)$ of an algebra over the Barratt-Eccles operad $A$
is endowed with Steenrod operations $\Sq^i: H^n(A)\,\rightarrow\,H^{n+i}(A)$
defined for $i\in\Z$.
If a class $c\in H^n(A)$ is represented by an element $a\in A^n$,
then the reduced square $\Sq^i(c)\in H^{n+i}(A)$
is represented by the product $\theta_{n-i}(a,a)\in A^{n+i}$.
The reduced squares $\Sq^i(c)$ where $i>n$ are $0$.
Let us observe that the top operation $\Sq^n(c)$ is identified with the Frobenius $\Sq^n(c) = c^2$.
In general, the Steenrod operations are quadratic in regard to scalar multiplication.
Explicitly, for $\lambda\in\F$, we have $\Sq^i(\lambda\,c) = \lambda^2\Sq^i(c)$.
On the other hand, for $c_1,c_2\in H^*(A)$, we have $\Sq^i(c_1+c_2) = \Sq^i(c_1)+\Sq^i(c_2)$.
In the case $A = N^*(X)$ and $H^*(A) = H^*(X,\F) = H^*(X,\F_2)\otimes_{\F_2}\F$,
we recover the classical Steenrod operations.
In this context, we have $\Sq^0(c\otimes\lambda)=c\otimes\lambda^2$
and $\Sq^i(c\otimes\lambda)=0$ for $i<0$.

In the general framework of $\Z$-graded Steenrod operations,
the Adem relation reads
$$\Sq^i\Sq^j = \sum\nolimits_{k\in\Z} \pmatrix{j-k-1\cr i-2k\cr}\Sq^{i+j-k}\Sq^k$$
(where $i<2j$ are given).
In fact, we may assume that the sum ranges over $i-j+1\leq k\leq i/2$
because the binomial coefficient vanishes outside this range.
Similarly, the Cartan formula reads
$$\Sq^n(a\cdot b) = \sum\nolimits_{k\in\Z} \Sq^k(a)\cdot\Sq^{n-k}(b).$$
This sum is also finite, because we have $\Sq^k(a) = 0$ for $k>\deg(a)$
and $\Sq^{n-k}(b) = 0$ for $n-k>\deg(b)$.

\tit{\it The big Steenrod algebra}\endtit
The big Steenrod algebra $B$ is generated by the operations $\Sq^i$, where $i\in\Z$.
We have explicitly $B = \F\langle\,\Sq^i,\ i\in\Z\,\rangle/(\hbox{Adem})$,
where $(\hbox{Adem})$ denotes the two-sided ideal generated by the Adem relations above.
The product of Steenrod monomials is extended to $\F$ coefficients
according to the relation
$\Sq^i\lambda = \lambda^2\Sq^i$ (for any $\lambda\in\F$).
The classical Steenrod algebra $A$
can be identified
with the quotient of the big Steenrod algebra $B$
by the two-sided ideal
generated by the relations $\Sq^0 = 1$ and $\Sq^i = 0$ for $i<0$.
In fact, according to Mandell,
the left ideal generated by $1+\Sq^0$ is a two-sided ideal in $B$
and contains the elements $\Sq^i$ such that $i<0$.
Consequently, we have $A = B/B(1+\Sq^0$).

\tit{\it Unstable modules}\label{UnstableModules}\endtit
The classical definition of an unstable module ({\it cf}. [\bibref{Sch}, Chapter 1])
is generalized in the context of modules over the big Steenrod algebra.
To be precise,
an unstable $B$-module (an unstable module over the big Steenrod algebra)
is a graded $\F$-module $M^*$
equipped with quadratic operations $\Sq^i: M^n\,\rightarrow\,M^{n+i}$,
(so that $M^*$ is a module over the algebra $B$)
such that $\Sq^i: M^n\,\rightarrow\,M^{n+i}$ vanishes for $i>n$.
Accordingly,
an unstable module over the big Steenrod algebra is equipped with the structure of a restricted $\F$-module.
By definition,
a restricted $\F$-module is a graded $\F$-module $M^*$ equipped with a Frobenius $\phi: M^*\,\rightarrow\,M^{2 *}$.
In the case of unstable modules,
we consider the top Steenrod operation $\Sq^n: M^n\,\rightarrow\,M^{2 n}$.

The free unstable $B$-module (respectively, the free unstable $A$-module)
generated by a graded $\F$-module $V$
is denoted either by $B\cdot V$ (respectively, $A\cdot V$)
or more simply by $B V$ (respectively, $A V$).
The free unstable module generated by a single element $e^n$ of degree $n$
is also denoted by $B e^n$
(respectively, $A e^n$).
According to Mandell,
we have a short exact sequence of left $B$-modules
$$0\,\rto{}{}\,B e^n\,\rto{1+\Sq^0}{}\,B e^n\,\rto{}{}\,A e^n\,\rto{}{}\,0$$
and, furthermore, this short exact sequence is split in the category of restricted modules.
We observe that this property holds for $n\in\Z$.
If $n<0$,
then we have $A e^n = 0$.
Hence,
in this case,
we just observe that the map $e^n\mapsto (1+\Sq^0)e^n$
gives an isomorphism of unstable $B$-modules.

\tit{\it Unstable algebras}\endtit
The cohomology of an algebra over the Barratt-Eccles operad
has the structure of an unstable $B$-algebra.
To be precise,
an unstable $B$-algebra
is an associative and commutative $\F$-algebra $R^*$
equipped with quadratic operations $\Sq^i: R^n\,\rightarrow\,R^{n+i}$
(so that $R^*$ is a module over the algebra $B$)
such that $\Sq^i: R^n\,\rightarrow\,R^{n+i}$ vanishes for $i>n$
and such that $\Sq^n: R^n\,\rightarrow\,R^{2 n}$ is identified with the Frobenius of $R$.
We assume in addition that the operations $\Sq^n$
verify the Cartan formula
with respect to the product of $R$.

An unstable $B$-module $M$ has an associated free unstable $B$-algebra $U(M)$.
This algebra has the same construction as the classical free unstable $A$-algebra
({\it cf}. [\bibref{Sch}, Section 3.8]).
Explicitly, the algebra $U(M)$ is the quotient of the symmetric algebra generated by $M$
by the ideal generated by the relations $\Sq^n(x) = x^2$,
where $x\in M^n$.
(Hence, the map $M\mapsto U(M)$ is given by a functor on the category of restricted $\F$-modules.)

The cohomology of a free $\E$-algebra is a free unstable $B$-algebra.
To be more explicit, for any differential graded $\F$-module $V$,
the canonical morphism $H^*(V)\,\rightarrow\,H^*(\E(V))$
induces an isomorphism
$$U(B H^*(V))\,\rto{\simeq}{}\,H^*(\E(V)).$$
Let us remind the reader that the cohomology of an Eilenberg-MacLane space
is a free unstable $A$-algebra.
Explicitly,
we have a fundamental class $e^n\in H^n(K(\Z/2,n),\F)$
together with an isomorphism
$$U(A e^n)\,\rto{\simeq}{}\,H^*(K(\Z/2,n),\F)$$
({\it cf}. [\bibref{Sch}, Section 9.8]).

\subsection{Models of Eilenberg-MacLane spaces and division functors}

\smallskip
We recall the construction of the cell algebra $F_n$.
As above,
we assume that the ground field $\F$ has characteristic $p= 2$.

\tit{\it Construction of the cofibrant model of an Eilenberg-MacLane space}\label{ModelEilenbergMacLane}\endtit
The cell algebra $F_n$ introduced by Mandell
is given by a cell extension of the form
$$\xymatrix{ \E(\F\,e^n)\ar[r]\ar[d]_{\phi_f} & \E(\F\,e^n\oplus\F\,b^{n-1})\ar[d] \\
\E(\F\,e^n)\ar[r] & F_n \\ }$$
The morphism $\phi_f: \E(\F\,e^n)\,\rightarrow\,\E(\F\,e^n)$
maps the element $e^n\in\E(\F\,e^n)$
to a representative of the class $e^n + \Sq^0(e^n)\in H^*(\E(\F e^n))$.
For instance, we have $\phi_f(e^n) = e^n + \theta_n(e^n,e^n)$.

This cell algebra $F_n$ is equivalent to an almost free algebra
such that $F_n = \E(\F\,e^n\oplus\F\,b^{n-1})$.
Furthermore,
we have a morphism
$F_n\,\rightarrow\,N^*(K(\Z/2,n))$,
which maps the element $e^n\in F_n$
to the fundamental class of $K(\Z/2,n)$.

\tit{\it The comparison arguments}\label{ComparisonArguments}\endtit
According to Mandell, the morphism $F_n\,\rightarrow\,N^*(K(\Z/2,n))$ is a quasi-iso\-mor\-phism
({\it cf}. [\bibref{Mnd}, Theorem 6.2]).
We outline Mandell's arguments.
We should observe that the result is valid for $n\in\Z$.
If $n<0$,
then we have $K(\Z/2,n) = \pt$ and we obtain simply $F_n\sim\F$.

We have $H^*(\E(\F\,e^n)) = U(B e^n)$.
We consider the morphism of free unstable algebras
$\psi_f: U(B e^n)\,\rightarrow\,U(B e^n)$
determined by $\phi_f: \E(\F\,e^n)\,\rightarrow\,\E(\F\,e^n)$.
One observes precisely that $\psi_f: U(B e^n)\,\rightarrow\,U(B e^n)$
is induced by the morphism of unstable $B$-modules
$f: B e^n\,\rightarrow\,B e^n$
such that $f(e^n) = e^n + \Sq^0(e^n)$.
We mention in paragraph \ref{UnstableModules}
that this morphism is split-injective in the category of restricted $\F$-modules.
More precisely,
if $f^*(B e^n)$ denotes the module obtained by restriction of structure,
then we have $f^*(B e^n) = B e^n\oplus A e^n$.
Consequently,
the morphism $\psi_f: U(B e^n)\,\rightarrow\,U(B e^n)$ makes $U(B e^n)$
a free module over $U(B e^n)$.
We have explicitly $\psi_f^*(U(B e^n)) = U(B e^n)\otimes U(A e^n)$.

The cocartesian square of paragraph \ref{ModelEilenbergMacLane}
gives rise to a spectral sequence $E^r\,\Rightarrow\,H^*(F_n)$
such that $E^2 = \Tor^{U(B e^n)}_*(\psi_f^*(U(B e^n)),\F)$.
We deduce from the result above that this spectral sequence degenerates,
because we obtain
$$\Tor^{U(B e^n)}_*(\psi_f^*(U(B e^n)),\F)
= \left\{\matrix{ U(A e^n),\hfill & \ \hbox{for}\ * = 0,\hfill\cr
0,\hfill & \ \hbox{otherwise}.\hfill\cr }\right.$$
Consequently,
the cohomology of $F_n$ is a free unstable $A$-algebra
and is generated by the class of the element $e^n\in F_n$.
One concludes readily
that the cohomology morphism $H^*(F_n)\,\rightarrow\,H^*(K(\Z/2,n))$
is an isomorphism.

\smallskip
We consider the image of the cell algebra $F_n$
under a division functor $-\oslash K: \E\Alg\,\rightarrow\,\E\Alg$,
where $K$ is a given coalgebra.
By lemma \ref{DivisionFreeMorphism},
we have:

\lem\label{DivisionEilenbergMacLane}\endtit{\it The algebra $F_n\oslash K$ fits in a cell extension
$$\xymatrix{ \E(\F\,e^n\otimes K)\ar[r]\ar[d]_{\phi_f\oslash K} & \E(\F\,e^n\otimes K\oplus\F\,b^{n-1}\otimes K)\ar[d] \\
\E(\F\,e^n\otimes K)\ar[r] & F_n\oslash K \\ }$$
Furthermore,
the morphism $f\oslash K: \F\,e^n\otimes K\,\rightarrow\,\E(\F\,e^n\otimes K)$
satisfies the relation
$$(f\oslash K)(e^n\otimes c) = e^n\otimes c
+ \sum_{k=0}^n (\tau^k\cdot\theta_{n-k})(e^n\otimes c^i_{(1)},e^n\otimes c^i_{(2)}),$$
for all $c\in K$, where $\theta_k{}^*(c) = \sum_i c^i_{(1)}\otimes c^i_{(2)}\in K^{\otimes 2}$.}

\tit{\it Construction of the comparison morphism}\label{ComparisonMorphismDivision}\endtit
We assume that $K$ is the chain complex of a simplicial set $K = N_*(X)$.
Furthermore,
we restrict ourself to the case of a finite simplicial set $X$.
We extend our constructions to infinite simplicial sets by limit arguments
({\it cf}. paragraph \ref{ComparisonMorphism}).

We construct a morphism
$F_n\oslash N_*(X)\,\rightarrow\,N_*\Map(X,K(\Z/2,n))$
natural in $X$.
We consider the universal resolution of $N^*(K(\Z/2,n))$
and the morphism
$F_{K(\Z/2,n)}\oslash N_*(X)\,\rightarrow\,N^*\Map(X,K(\Z/2,n))$
introduced in paragraph \ref{ComparisonMorphism}.
We have a quasi-iso\-mor\-phism $F_n\,\rightarrow\,F_{K(\Z/2,n)}$,
because $F_n$ and $F_{K(\Z/2,n)}$
are both cofibrant resolutions of $N^*(K(\Z/2,n))$.
Therefore,
we have a natural comparison morphism $F_n\oslash N_*(X)\,\rightarrow\,N_*\Map(X,K(\Z/2,n))$,
which is given by the composite
$$F_n\oslash N_*(X)\,\rightarrow\,F_{K(\Z/2,n)}\oslash N_*(X)\,\rightarrow\,N^*\Map(X,K(\Z/2,n)).$$

We generalize the arguments outlined in paragraph \ref{ComparisonArguments}
in order to prove that this morphism is a quasi-isomorphism.
We still assume that the ground field has characteristic $p = 2$.
In fact,
we have fixed representatives of the Steenrod squares $\Sq^n$, $n\in\Z$,
but we just need to assume the Cartan relation
$\Sq^0(a\cdot b) = \sum_{n\in\Z} \Sq^{n}(a)\cdot\Sq^{-n}(b)$.
Therefore,
the generalization of our arguments to odd primes is straighforward.

To begin with,
we obtain the following result:

\lem\label{CohomologyDivision}\endtit
{\it We have $H^*\bigl(\E(\F\,e^n\otimes N_*(X))\bigr) = U\bigl(B(\F\,e^n\otimes H_*(X))\bigr)$.
The morphism of free unstable algebras
$$\displaylines{\psi_g: U\bigl(B(\F\,e^n\otimes H_*(X))\bigr)
\,\rightarrow\,U\bigl(B(\F\,e^n\otimes H_*(X))\bigr) \cr
\noalign{\hbox{determined by}}
\phi_f\oslash N_*(X): \E(\F\,e^n\otimes N_*(X))\,\rightarrow\,\E(\F\,e^n\otimes N_*(X)) \cr
\noalign{\hbox{is induced by the morphism of unstable $B$-modules}}
g: B(\F\,e^n\otimes H_*(X))\,\rightarrow\,B(\F\,e^n\otimes H_*(X)) \cr
\noalign{\hbox{such that}}
g(e^n\otimes c) = e^n\otimes c + \sum\nolimits_{l\in\Z} \Sq^l(e^n\otimes(c\Sq^{-l})), \cr }$$
for all $c\in H_*(X)$.}

\proof\endtit
We would like to determine the cohomology of the morphism
$$f': \F\,e^n\otimes N_*(X)\,\rightarrow\,\bigl(\E(2)\otimes(\F\,e^n\otimes N_* X)^{\otimes 2}\bigr)_{\Sigma_2}$$
such that
$f'(e^n\otimes x) = \sum_{k=0}^n (\tau^k\cdot\theta_{n-k})(e^n\otimes x^i_{(1)},e^n\otimes x^i_{(2)})$,
where $\theta_k^*(x) = \sum_i x^i_{(1)}\otimes x^i_{(2)}$
({\it cf}. lemma \ref{DivisionEilenbergMacLane}).
We consider homological Steenrod operations $\Sq^k: H_m(X)\,\rightarrow\,H_{m-k}(X)$.
These operations are deduced from the cohomological ones by a duality process.
Therefore,
we introduce the dual
$$f'{}^{\vee}: \bigl(\E(2)^{\vee}\otimes(\F\,e_n\otimes N^* X)^{\otimes 2}\bigr)^{\Sigma_2}\,\rightarrow\,\F\,e_n\otimes N^*(X)$$
of the morphism $f'$.
We let $(c_i)_i$ denote a basis of $H_*(X,\F_2)$.
We fix a representative $x_i\in N_*(X)$ of each class $c_i\in H_*(X)$.
We consider dual basis elements $c^i\in H^*(X)$ and dual cocycles $x^i\in N^*(X)$.
We have also $\E(2)^{\vee} = \F[\Sigma_2]\,\theta^k$,
where $\theta^k\in\E(2)^{\vee}$ is dual to $\theta_k\in\E(2)$
and satisfies the relation
$\delta(\theta^k) = \theta^{k+1}+\tau\theta^{k+1}$.

The cohomology module
$$H^*\bigl(\E(2)\otimes(e^n\otimes N_* X)^{\otimes 2}\bigr)_{\Sigma_2}
\subset H^*\bigl(\E(\F\,e^n\otimes N_*(X))\bigr) = U\bigl(B(\F\,e^n\otimes H^*(X))\bigr)$$
is generated by the cycles $\theta_0(e^n\otimes x_i,e^n\otimes x_j)$
and $\theta_{n-k}(e^n\otimes x_i,e^n\otimes x_i)$,
where $k\leq n$.
The former represent the products $e^n\otimes c_i\cdot e^n\otimes c_j\in U\bigl(B(\F\,e^n\otimes H^*(X))\bigr)$
and the latter represent the Steenrod squares $\Sq^l(e^n\otimes c_i)\in U\bigl(B(\F\,e^n\otimes H^*(X))\bigr)$,
where $l = k - \deg(c_i)\leq n - \deg(c_i)$.
Dually,
the homology module
$$H_*\bigl(\E(2)^{\vee}\otimes(\F\,e_n\otimes N^* X)^{\otimes 2}\bigr)^{\Sigma_2}
\subset H_*\bigl(\E(\F\,e^n\otimes N_*(X))^{\vee}\bigr)$$
is generated by the cycles
$(1+\tau)\theta^0(e_n\otimes x^i,e_n\otimes x^j) + (1+\tau)\theta^0(e_n\otimes x^j,e_n\otimes x^i)$
and $(1+\tau)\theta^{n-k}(e_n\otimes x^i,e_n\otimes x^i)$,
where $k\leq n$.

We obtain
$$\displaylines{ f'{}^{\vee}\bigl((1+\tau)\theta^0(e_n\otimes x^i,e_n\otimes x^j)
+ (1+\tau)\theta^0(e_n\otimes x^j,e_n\otimes x^i)\bigr)\hfill\cr
\hfill = e_n\otimes\theta_n(x^i,x^j) + e_n\otimes\theta_n(x^j,x^i). \cr }$$
The element $\theta_n(x^i,x^j) + \theta_n(x^j,x^i)\in N^*(X)$ vanishes in $H^*(X)$,
because it corresponds to the differential of $\theta_{n+1}(x^i,x^j)\in N^*(X)$.
We obtain also
$$f'{}^{\vee}\bigl((1+\tau)\theta^{n-k}(e_n\otimes x^i,e_n\otimes x^i)\bigr)
= \left\{\matrix{ e_n\otimes\theta_k(x^i,x^i),\hfill & \ \hbox{for}\ 0\leq k\leq n,\hfill\cr
0,\hfill & \ \hbox{for}\ k< 0,\hfill \cr }\right.$$
The element $\theta_k(x^i,x^i)\in N^*(X)$
is a representative of the Steenrod square $\Sq^{-l}(c^i)\in H^*(X)$,
where $l = k - \deg(c^i)$,
as in the paragraph above.
We have $\Sq^{-l}(c^i) = \sum_j \lambda^{-l}_{i j} c^j$,
where $\deg(c^j) = 2\deg(c^i) - k = - 2l + k$.

By duality,
the homology morphism $H_*(f')$
maps the element
$e^n\otimes c_j\in\F\,e^n\otimes H_*(X)$
to the sum
$\sum_l \lambda^{-l}_{i j} \Sq^l(e^n\otimes c_i)\in H^*\bigl(\E(\F\,e^n\otimes N_*(X))\bigr)$.
By definition of homological Steenrod operations,
we have $c_j\Sq^{-l} = \sum_i \lambda^{-l}_{i j} c_i$.
Accordingly,
the sum above is an expansion of the expression $\sum_l \Sq^l(e^n\otimes c_j\Sq^{-l})$.
The summation index $l$ is such that $\deg(c_j) = - 2l + k$.
Since $0\leq k\leq n$, we obtain $(-\deg(c_j))/2\leq l\leq(n-\deg(c_j))/2$.
But,
the unstability relations implies that the formula $\Sq^l(e^n\otimes(c_j\Sq^{-l}))$
vanishes outsider this range.
Therefore,
we can assume that the sum ranges over $\Z$
and we obtain the result stated in the lemma.

\lem\label{RestrictedSplitting}\endtit
{\it The morphism $g: B(\F\,e^n\otimes H_*(X))\,\rightarrow\,B(\F\,e^n\otimes H_*(X))$
is split-injective in the category of restricted $\F$-modules
and fits in a short exact sequence of unstable $B$-modules
$0\,\rto{}{}\,B(\F\,e^n\otimes H_*(X))\,\rto{}{}\,B(\F\,e^n\otimes H_*(X))
\,\rto{}{}\,A(\F\,e^n\otimes H_*(X))\,\rto{}{}\,0$.}

\proof\endtit
The sequence
$$0\,\rto{}{}\,B(e^n\otimes c)\,\rto{1+\Sq^0}{}\,B(e^n\otimes c)\,\rto{}{}\,A(e^n\otimes c)\,\rto{}{}\,0,$$
is exact,
for all $c\in H_*(X)$ ({\it cf}. paragraph \ref{UnstableModules}).
Furthermore,
we have a morphism of restricted modules $r_0: B(e^n\otimes c)\,\rightarrow\,B(e^n\otimes c)$
such that $r_0(b(1+\Sq^0)(e^n\otimes c)) = b(e^n\otimes c)$,
for all $b(e^n\otimes c)\in B(\F\,e^n\otimes H_*(X))$.
To prove our property,
we observe that the morphism
$g: B(\F\,e^n\otimes H_*(X))\,\rightarrow\,B(\F\,e^n\otimes H_*(X))$
is represented by a triangular matrix
which has the element $(1+\Sq^0)$
on all diagonal entries.

To be more explicit,
we construct a morphism of restricted $\F$-modules
$r: B(\F\,e^n\otimes H_*(X))\,\rightarrow\,B(\F\,e^n\otimes H_*(X))$
such that $r\,g(b(e^n\otimes c)) = b(e^n\otimes c)$
by induction on the degree of $c\in H_*(X)$.
Since $c\Sq^{-l} = 0$ for $l>0$ and $c\Sq^0 = c$,
we obtain
$$\eqalign{ g(b(e^n\otimes c)) & = b(1+\Sq^0)(e^n\otimes c) + g'(b(e^n\otimes c)), \cr
\hbox{where}\qquad g'(b(e^n\otimes c)) & = \sum\nolimits_{l<0} b\Sq^l(e^n\otimes c\Sq^{-l}). \cr }$$
In this equation,
we have $\deg(c\Sq^{-l}) = \deg(c) + l < \deg(c)$.
Therefore, by induction, we can set
$r(b(e^n\otimes c)) = r_0(b(e^n\otimes c)) - r\,g'(r_0(b(e^n\otimes c))$.
We obtain readily the identity $r\,g(b(e^n\otimes c)) = b(e^n\otimes c)$.

Similarly,
a straightforward induction proves that the cokernel
of the morphism $g: B(\F\,e^n\otimes H_*(X))\,\rightarrow\,B(\F\,e^n\otimes H_*(X))$
is the module $A(\F\,e^n\otimes H_*(X))$.

\smallskip
As in Mandell's proof ({\it cf}. paragraph \ref{ComparisonArguments}),
lemma \ref{RestrictedSplitting} has the following formal consequence:

\lem\endtit{\it Let $U_B^X = U\bigl(B(\F\,e^n\otimes H_*(X))\bigr)$
and $U_A^X = U\bigl(A(\F\,e^n\otimes H_*(X))\bigr)$.
The morphism $\psi_g: U_B^X\,\rightarrow\,U_B^X$ makes the algebra $U_B^X$ a free module over $U_B^X$.
We have $\psi_g^* U_B^X = U_B^X\otimes U_A^X$.}

\tit{\it The spectral sequence}\endtit
The cell extension of lemma \ref{DivisionEilenbergMacLane}
gives rise to a spectral sequence $E^r\,\Rightarrow\,H^*(F_n\oslash N_*(X))$
such that $E^2 = \Tor_*^{U_B^X}(\psi_g^* U_B^X,\F)$.
The lemma above implies that this spectral sequence degenerates at $E^2$,
because we obtain
$$\Tor^{U_B^X}_*(\psi_f^* U_B^X,\F)
= \left\{\matrix{ U_A^X,\hfill & \ \hbox{for}\ * = 0,\hfill\cr
0,\hfill & \ \hbox{otherwise}.\hfill\cr }\right.$$
Consequently:

\lem\endtit{\it We have an isomorphism
$U\bigl(A(\F\,e^n\otimes H_*(X))\bigr)\,\rto{\simeq}{}\,H^*(F_n\oslash N_*(X))$.}

\smallskip
We can achieve the proof of the main result of this section:

\lem\endtit{\it The comparison morphism $F_n\oslash N_*(X)\,\rightarrow\,N^*\Map(X,K(\Z/2,n))$
(defined in paragraph \ref{ComparisonMorphismDivision}) is a quasi-isomorphism.}

\proof\endtit
We have also an isomorphism
$U\bigl(A(\F\,e^n\otimes H_*(X))\bigr)\,\rto{\simeq}{}\,H^*\bigl(\Map(X,K(\Z/2,n))\bigr)$,
because $\Map(X,K(\Z/2,n))$ is a generalized Eilenberg-MacLane space.
We observe that the comparison morphism preserves the cohomology generators.
For that purpose,
we consider the adjoint morphism $F_n\,\rightarrow\,N^*(X)\otimes N^*\Map(X,K(\Z/2,n))$.
According to the definition,
this morphism occurs as a composite
in a diagram
$$\xymatrix{ F_n\ar[r]\ar[dr] & F_Y\ar[d]_{\sim}\ar[r] & F_{X\times\Map(X,Y)}\ar[d]_{\sim}\ar[rd] & \\
& N^*(Y)\ar[r] & N^*(X\times\Map(X,Y))\ar[r] & N^*(X)\otimes N^*(\Map(X,Y)) \\ }$$
where the right triangle commutes in cohomology ({\it cf}. theorem \ref{EZKunnethEquivalence}).
Consequently,
the element $e^n\in F_n$
is mapped to a representative of the image of the fundamental class of $K(\Z/2,n)$
under the classical morphism
$$H^*(K(\Z/2,n))\,\rightarrow\,H^*(X)\otimes H^*\Map(X,K(\Z/2,n))$$
induced by the evaluation product.
By adjunction,
we deduce that the comparison morphism $F_n\oslash N_*(X)\,\rightarrow\,N^*\Map(X,K(\Z/2,n))$
identifies a cycle $e^n\otimes x\in e^n\otimes N_*(X)$
with a generator of $H^*\Map(X,K(\Z/2,n))$.
The conclusion follows.

\subsection{Fibrations and division functors}

\smallskip
The next lemma allows to achieve the proof of theorem \ref{resolutionMapSpace} by a classical induction process
({\it cf}. [\bibref{FS}], [\bibref{Sch}, Section 9.9]).

\lem\label{FibrationDivision}\endtit
{\it We are given a fibration sequence $F\,\rightarrow\,E\,\rightarrow\,B$
where $B$, $E$, $F$ are $p$-$\pi_*$-finite ({\it cf}. [\bibref{Sch}, Section 9.7]).
If the comparison morphism $F_Y\oslash N_*(X)\,\rightarrow\,N^*(\Map(X,Y))$
is a quasi-isomorphism for $Y = B$ and $Y = E$,
then the comparison morphism is also a quasi-isomorphism for $Y = F$.}

\proof\endtit
We consider the natural morphism $F_B\,\rightarrow\,F_E$,
which is a cofibration by proposition \ref{UniversalCellMorphism}.
We form the cocartesian square
$$\xymatrix{ F_B\ar[d]\ar[r] & F_E\ar[d] \\
\F\ar[r] & \F\vee_{F_B} F_E \\ }$$
We have also a morphism $F_E\,\rightarrow\,F_F$
by functoriality of universal cofibrant resolutions.
We claim that the induced morphism $\F\vee_{F_B} F_E\,\rightarrow\,F_F$
is a quasi-isomorphism.
This property follows from a result of Mandell ({\it cf}. [\bibref{Mnd}, Lemma 5.2]).
Explicitly,
the algebra $\F\vee_{F_B} F_E$ is quasi-isomorphic to $N^*(F)$.
Hence,
the pushout process above provides a cofibrant resolution of $N^*(F)$,
which is necessarily equivalent to the universal one $F_F$.
The division functor $-\oslash N_*(X): \E\Alg\,\rightarrow\,\E\Alg$
preserves cocartesian squares
and quasi-isomorphisms of cofibrant algebras,
such as $\F\vee_{F_B} F_E\,\rightarrow\,F_F$.
Therefore,
we obtain a homotopy cocartesian square
$$\xymatrix{ F_B\oslash N_*(X)\ar[d]\ar[r]\ar@{}[dr]|{(M)} & F_E\oslash N_*(X)\ar[d] \\
\F\ar[r] & F_F\oslash N_*(X) \\ }$$
This homotopy cocartesian square gives rise to a strongly convergent spectral sequence
$E^r_{(M)}\,\rightarrow\,H^*(F_F\oslash N_*(X))$
such that $E^2_{(M)} = \Tor^{H^*(F_B\oslash N_*(X))}(H^*(F_E\oslash N_*(X)),\F)$
({\it cf}. [\bibref{Mnd}, Corollary 3.6]).

We compare this spectral sequence to the Eilenberg-Moore spectral sequence
$E^r_{(C)}\,\rightarrow\,H^*(\Map(X,F))$
of the fibration
$\Map(X,F)\,\rightarrow\,\Map(X,E)\,\rightarrow\,\Map(X,B)$.
We have $E^2_{(C)} = \Tor^{H^*(\Map(X,B))}(H^*(\Map(X,E)),\F)$
and the finiteness assumption of the lemma ensures the strong convergence of the spectral sequence.
By functoriality,
the comparison morphism gives a morphism
from the commutative square of cofibrant models $(M)$
to the commutative square of cochain algebras $(C)$
$$\xymatrix{ N^*(\Map(X,B))\ar[d]\ar[r]\ar@{}[dr]|{(C)} & N^*(\Map(X,E))\ar[d] \\
\F\ar[r] & N^*(\Map(X,F)) \\ }$$
As a consequence,
we have a morphism of spectral sequence $E^r_{(M)}\,\rightarrow\,E^r_{(C)}$.
If we have isomorphisms
$$H^*(F_B\oslash N_*(X))\,\rightarrow\,H^*(\Map(X,B))
\qquad\hbox{and}
\qquad H^*(F_E\oslash N_*(X))\,\rightarrow\,H^*(\Map(X,E)),$$
then the $E^2$-stages of the spectral sequences agree.
We conclude that the natural morphism
$H^*(F_F\oslash N_*(X))\,\rightarrow\,H^*(\Map(X,F))$
is an isomorphism.
This achieves the proof of lemma \ref{FibrationDivision}.

\section{The Eilenberg-Zilber equivalence}

\smallskip
We consider the classical shuffle morphism
$N^*(X\times Y)\,\rightarrow\,N^*(X)\widehat{\otimes} N^*(Y)$
({\it cf}. [\bibref{MLH}, Section VIII.8]),
where $N^*(X)\widehat{\otimes} N^*(Y)$ is the profinite completion of the tensor product $N^*(X)\otimes N^*(Y)$.
We observe that this morphism of dg-modules
is not a morphism of $\E$-algebras.
Therefore,
the purpose of this section is to construct a homotopy Eilenberg-Zilber morphism
from the universal almost free algebra $F_{X\times Y}$
to the complete tensor product $N^*(X)\widehat{\otimes} N^*(Y)$.
We obtain the following result:

\th\label{EZKunnethEquivalence}\endtit{\it Let $X$ and $Y$ be simplicial sets.
We have a morphism of $\E$-algebras
$F_{X\times Y}\,\rightarrow\,N^*(X)\widehat{\otimes} N^*(Y)$,
functorial in $X$ and $Y$,
whose cohomology corresponds to the classical K\"unneth isomorphism
$H^*(X\times Y)\rto{\simeq}{}\,H^*(X)\widehat{\otimes} H^*(Y)$.
Precisely,
the diagram
$$\xymatrix{ F_{X\times Y}\ar[dr]\ar[d]_{\sim} & \\
N^*(X\times Y)\ar[r]^{\sim} & N^*(X)\widehat{\otimes} N^*(Y) \\ }$$
where $N^*(X\times Y)\,\rto{\sim}{}\,N^*(X)\widehat{\otimes} N^*(Y)$ is the classical shuffle morphism,
gives rise to a commutative diagram in cohomology.}

\smallskip
We have $N^*(X)\widehat{\otimes} N^*(Y) = \lim_{K,L} N^*(K)\widehat{\otimes} N^*(L)$,
where $K$ (respectively, $L$) ranges over finite simplicial subsets of $X$ (respectively, $Y$).
Therefore,
we reduce to the case of finite simplicial sets
and omit profinite completions.
We construct a homotopy Eilenberg-Zilber morphism
as stated in lemma \ref{EZDefinition} below.
We prove that the assertion
of theorem \ref{EZKunnethEquivalence}
about the K\"unneth isomorphism
follows from the construction
({\it cf}. lemma \ref{EZKunneth}).

\lem\label{EZDefinition}\endtit{\it Let $X$ and $Y$ be finite simplicial sets.
There is a morphism of $\E$-algebras
$$F_{X\times Y}\,\rightarrow\,N^*(X)\otimes N^*(Y)$$
which is functorial in $X$ and $Y$
and which reduces to the augmentation morphism $F_\pt\,\rightarrow\,N^*(\pt)$
for $X = Y = \pt$.}

\tit{\it Reminder: the structure of the universal almost free resolution}\endtit
We recall that the universal almost free resolution of an algebra $A$
is given by the construction $F_A = \E(C_A)$,
where $C_A = \B\E(A)$ is an almost cofree coalgebra over the bar cooperad $\B\E$.
This coalgebra is equipped with a differential $\delta_A: \B\E(A)\,\rightarrow\,\B\E(A)$
which differs from the canonical differential of the cofree coalgebra $\delta: \B\E(A)\,\rightarrow\,\B\E(A)$
by a coderivation $d_A: \B\E(A)\,\rightarrow\,\B\E(A)$
determined by the algebra structure of $A$.
We have $\B\E(A) = \bigoplus_r \B\E_{(r)}(A)$,
where $\B\E_{(r)}(A) = (\B\E(r)\otimes A^{\otimes r})_{\Sigma_r}$.
We have a natural isomorphism $A\simeq\B\E_{(1)}(A)$
and the canonical surjection $\B\E(A)\,\rightarrow\,\B\E_{(1)}(A) = A$
is identified with the universal morphism $\B\E(A)\,\rightarrow\,A$.
If we consider the canonical section $A = \B\E_{(1)}(A)\,\rightarrow\,\B\E(A)$,
then we obtain a morphism of dg-modules $A\,\rightarrow\,C_A$
which identifies the algebra $A$
with the indecomposable part of the almost cofree coalgebra $C_A$.

By definition of the universal almost free resolution,
a morphism of $\E$-algebras $F_A\,\rightarrow\,B$
is equivalent to a morphism of $\B\E$-coalgebras $C_A\,\rightarrow\,C_B$.
Furthermore,
the canonical morphism $F_A\,\rightarrow\,A$
corresponds to the identity of $C_A$.
In particular,
the morphism of lemma \ref{EZDefinition}
is equivalent to a morphism of $\B\E$-coalgebras
$C_{N^*(X\times Y)}\,\rightarrow\,C_{N^*(X)\otimes N^*(Y)}$,
which reduces to the identity morphism for $X = Y = \pt$.

\smallskip
We would like to mention that the morphism $C_{N^*(X\times Y)}\,\rightarrow\,C_{N^*(X)\otimes N^*(Y)}$
is an extension of the classical Eilenberg-Zilber equivalence.
Precisely,
we obtain the following result:

\prop\label{EZExtension}\endtit{\it The morphism of $\B\E$-coalgebras
$C_{N^*(X\times Y)}\,\rightarrow\,C_{N^*(X)\otimes N^*(Y)}$,
deduced from lemma \ref{EZDefinition},
fits in a commutative diagram of dg-modules
$$\xymatrix{ N^*(X\times Y)\ar[d]\ar[r] & N^*(X)\otimes N^*(Y)\ar[d] \\
C_{N^*(X\times Y)}\ar[r] & C_{N^*(X)\otimes N^*(Y)} \\ }$$
where $N^*(X\times Y)\,\rightarrow\,N^*(X)\otimes N^*(Y)$ is the classical shuffle morphism.}

\proof\endtit
Recall that $N^*(X\times Y)$ (respectively, $N^*(X)\otimes N^*(Y)$)
is identified with the indecomposable part of $C_{N^*(X\times Y)}$
(respectively, $C_{N^*(X)\otimes N^*(Y)}$).
Because of this property,
the morphism of coalgebras $C_{N^*(X\times Y)}\,\rightarrow\,C_{N^*(X)\otimes N^*(Y)}$
induces a morphism of dg-modules $N^*(X\times Y)\,\rightarrow\,N^*(X)\otimes N^*(Y)$
as in the diagram above.
The induced morphism coincides with the classical shuffle morphism
because it verifies the same characteristic properties.
Namely,
we obtain a morphism $N^*(X\times Y)\,\rightarrow\,N^*(X)\otimes N^*(Y)$
which is functorial in $X$ and $Y$
and which reduces to the identity morphism
for $X = Y = \pt$
({\it cf}. [\bibref{DEZ}]).

\smallskip
We prove the assertion of theorem \ref{EZKunnethEquivalence}
about the K\"unneth isomorphism
by similar arguments:

\lem\label{EZKunneth}\endtit{\it The morphism provided by lemma \ref{EZDefinition}
makes the following diagram commutative in cohomology
$$\xymatrix{ F_{N^*(X\times Y)}\ar[dr]\ar[d]_{\sim} & \\
N^*(X\times Y)\ar[r]^{\sim} & N^*(X)\otimes N^*(Y) \\ }$$}

\proof\endtit
The composite of the canonical injections
$N^*(X\times Y)\,\rightarrow\,\B\E(N^*(X\times Y)) = C_{N^*(X\times Y)}$
and
$C_{N^*(X\times Y)}\,\rightarrow\,\E(C_{N^*(X\times Y)}) = F_{N^*(X\times Y)}$
defines a morphism of dg-modules $N^*(X\times Y)\,\rightarrow\,F_{N^*(X\times Y)}$.
This morphism is a section of the canonical augmentation $F_{N^*(X\times Y)}\,\rightarrow\,N^*(X\times Y)$.
Hence,
the induced morphism
$H^*(N^*(X\times Y))\,\rightarrow\,H^*(F_{N^*(X\times Y)})$
is an inverse
of the cohomology isomorphism $H^*(F_{N^*(X\times Y)})\,\rto{\simeq}{}\,H^*(N^*(X\times Y))$.
The lemma follows from the observation that the composite
$$N^*(X\times Y)\,\rightarrow\,F_{N^*(X\times Y)}\,\rightarrow\,N^*(X)\otimes N^*(Y)$$
coincides with the shuffle morphism.
As in the proof of lemma \ref{EZDefinition},
this assertion is an immediate consequence
of the characteristic properties
of the Eilenberg-Zilber equivalence.

\smallskip
We define the morphism
$F_{N^*(X\times Y)}\,\rightarrow\,N^*(X)\otimes N^*(Y)$
and prove lemma \ref{EZDefinition}
in the next paragraphs.
We adapt the classical methods of acyclic models to our framework.
First,
we restrict to the case $X = \Delta^m$ and $Y = \Delta^n$.
To be precise,
the full subcategory of $\Spaces^2$
formed by the pairs $(\Delta^m,\Delta^n)$, $(m,n)\in\N^2$
is isomorphic to $\Delta^2$.
Hence,
we consider the covariant functor $\Delta^2\,\rightarrow\,\Spaces^2$
provided by the map $(m,n)\mapsto(\Delta^m,\Delta^n)$.
We construct a functorial morphism of algebras
$F_{N^*(\Delta^m\times\Delta^n)}\,\rightarrow\,N^*(\Delta^m)\otimes N^*(\Delta^n)$
by induction on dimension.
For that purpose,
we introduce Reedy structures in the framework of dg-modules and $\E$-algebras
({\it cf}. [\bibref{R}]).
We refer to the monographs [\bibref{DHK}], [\bibref{GoJ}] and [\bibref{H}]
for more background about Reedy's theory.

\tit{\it Bisimplicial modules and Reedy structures}\endtit
Let $V$ be a bisimplicial (differential graded) module.
The matching module $M_{m,n} V$ is formed by sequences of elements
$v_i\in V_{m-1,n}$, $0\leq i\leq m$, and $w_j\in V_{m,n-1}$, $0\leq j\leq n$,
such that $d_k\times 1(v^i) = d_{i-1}\times 1(v_k)$, for $k<i$,
and $1\times d_l(w^j) = 1\times d_{j-1}(w^l)$, for $l<j$.
We have a canonical morphism $V_{m,n}\,\rightarrow\,M_{m,n} V$
which maps an element $x\in V_{m,n}$
to $(d_0\times 1(x),\ldots,d_m\times 1(x),1\times d_0(x),\ldots,d_n(x))\in M_{m,n} V$.
We have equivalently $M_{m,n} V = \lim_{(u,v)} V_{m',n'}$,
where the limit ranges over proper pairs of injective morphisms
$u: [m']\,\rightarrow\,[m]$ and $v: [n']\,\rightarrow\,[n]$
in the simplicial category.

Dually,
we have a latching module $L_{m,n} V$ together with a morphism $L_{m,n} V\,\rightarrow\,V_{m,n}$.
The latching module verifies $L_{m,n} V = \colim_{(u,v)} V_{m',n'}$,
where the colimit ranges over proper pairs of surjective morphisms
$u: [m]\,\rightarrow\,[m']$ and $v: [n]\,\rightarrow\,[n']$.

We have the following property:

\tit{\sc Fact}\label{ReedyModules}\endtit
{\it All bisimplicial dg-modules are Reedy fibrant and Reedy cofibrant.
More explicitly,
the matching morphism $V_{m,n}\,\rightarrow\,M_{m,n} V$ is always surjective
and
the latching morphism $L_{m,n} V\,\rightarrow\,V_{m,n}$ is always injective.}

\proof\endtit
The proof of these assertions is classical.
Suppose given $(v_0,\ldots,v_m,w_0,\ldots,w_n)\in M_{m,n} V$.
We construct an element $x\in V_{m,n}$
such that
$(d_0\times 1(x),\ldots,d_m\times 1(x),1\times d_0(x),\ldots,1\times d_n(x)) = (v_0,\ldots,v_m,w_0,\ldots,w_n)$.
We assume by induction $v_0 = \cdots = v_{i-1} = 0$.
We have an element $x'\in V_{m,n}$
such that $d_0\times 1(x') = \cdots = d_{i-1}\times 1(x') = 0$ and $d_i\times 1(x') = v_i$.
For instance,
we can prove by straighforward verifications
that the element
$x' = (-1)^i s_0\times 1(v_i) + (-1)^{i-1} s_1\times 1(v_i) + \cdots + s_i\times 1(v_i)$
satisfies these equations.
We set
$(v'_0,\ldots,v'_m,w'_0,\ldots,w'_n)
= (v_0,\ldots,v_m,w_0,\ldots,w_n)
- (d_0\times 1(x'),\ldots,d_m\times 1(x'),1\times d_0(x'),\ldots,1\times d_n(x'))$
and we resume the process.
Eventually,
we prove that the element $(v_0,\ldots,v_m,w_0,\ldots,w_n)\in M_{m,n} V$
lies in the image of the matching morphism.

Hence,
we conclude that the matching morphism $V_{m,n}\,\rightarrow\,M_{m,n} V$ is surjective.
By a similar construction,
we prove
that the matching morphism of a bicosimplicial module $V^{m,n}\,\rightarrow\,M^{m,n} V$ is surjective.
By duality,
this proves also
that the latching morphism of a bisimplicial module $L_{m,n} V\,\rightarrow\,V_{m,n}$ is injective.

\lem\label{ReedyAlgebras}\endtit
{\it The bisimplicial algebra $(m,n)\mapsto N^*(\Delta^m)\otimes N^*(\Delta^n)$ is Reedy fibrant.
Explicitly, given $(m,n)\in\N^2$,
the limit morphism
$$N^*(\Delta^m)\otimes N(\Delta^n)\,\rightarrow\,\lim_{(u,v)} N^*(\Delta^{m'})\otimes N^*(\Delta^{n'}),$$
where $u: [m']\,\rightarrow\,[m]$ and $v: [n']\,\rightarrow\,[n]$ range over proper pairs of injections,
is a fibration of $\E$-algebras.

The bisimplicial algebra $(m,n)\mapsto F_{N^*(\Delta^m\times\Delta^n)}$ is acyclic and Reedy cofibrant.
Explicitly, given $(m,n)\in\N^2$,
the colimit morphism
$$\colim_{(u,v)} F_{N^*(\Delta^{m'}\times\Delta^{n'})}\,\rightarrow\,F_{N^*(\Delta^m\times\Delta^n)},$$
where $u: [m]\,\rightarrow\,[m']$ and $v: [n]\,\rightarrow\,[n']$ range over proper pairs of surjections,
is an acyclic cofibration of $\E$-algebras.}

\proof\endtit
We deduce immediately from the case of dg-modules ({\it cf}. fact \ref{ReedyModules})
that the matching morphism
$$N^*(\Delta^m)\otimes N(\Delta^n)\,\rightarrow\,\lim_{(u,v)} N^*(\Delta^{m'})\otimes N^^*(\Delta^{n'}),$$
is a fibration,
because, in the category of $\E$-algebras, the fibrations are the surjective morphisms.

Consider the case of the bisimplicial algebra $(m,n)\mapsto F_{N^*(\Delta^m\times\Delta^n)}$.
The morphisms $\F\,\rightarrow\,F_{N^*(\Delta^m\times\Delta^n)}$ are clearly weak-equivalences.
Hence, we have just to prove that these morphisms form a Reedy cofibration.
We have
$$\colim_{(u,v)} F_{N^*(\Delta^{m'}\times\Delta^{n'})}
= \E(\colim_{(u,v)} C_{N^*(\Delta^{m'}\times\Delta^{n'})}),$$
by construction of colimits
in the context of almost free algebras.
Furthermore,
the latching morphism
$\colim_{(u,v)} F_{N^*(\Delta^{m'}\times\Delta^{n'})}\,\rightarrow\,F_{N^*(\Delta^m\times\Delta^n)}$
is induced by the morphism of dg-modules
$\colim_{(u,v)} C_{N^*(\Delta^{m'}\times\Delta^{n'})}\,\rightarrow\,C_{N^*(\Delta^m\times\Delta^n)}$.
We deduce from the case of dg-modules ({\it cf}. fact \ref{ReedyModules})
that this morphism is injective.
Our assertion follows from the results
of section \ref{CellFreeStructure}.
Namely,
since $F_{N^*(\Delta^m\times\Delta^n)} = \E(C_{N^*(\Delta^m\times\Delta^n)})$ is a cell algebra
({\it cf}. proposition \ref{UniversalCellResolution}),
lemma \ref{CellStrictMorphism} implies that the latching morphism
$\colim_{(u,v)} F_{N^*(\Delta^{m'}\times\Delta^{n'})}\,\rightarrow\,F_{N^*(\Delta^m\times\Delta^n)}$
is a relative cell inclusion
and therefore is a cofibration of $\E$-algebras.

\lem\label{EZModel}\endtit
{\it Suppose given compatible morphisms of algebras
$$F_{N^*(\Delta^{m'}\times\Delta^{n'})}\,\rightarrow\,N^*(\Delta^{m'})\otimes N^*(\Delta^{n'})$$
for all bisimplices $\Delta^{m'}\times\Delta^{n'}$,
such that $m'+n'<N$.
We consider a cartesian product $\Delta^m\times\Delta^n$,
where $m + n = N$.
There is a morphism
$$F_{N^*(\Delta^m\times\Delta^n)}\,\rightarrow\,N^*(\Delta^m)\otimes N^*(\Delta^n)$$
which makes all functoriality diagrams commute
for pairs of morphisms
$(u: [m']\,\rightarrow\,[m],v: [n']\,\rightarrow\,[n])$
and
$(u: [m]\,\rightarrow\,[m'],v: [n]\,\rightarrow\,[n'])$
such that $m'+n'<N$.}

\proof\endtit
The given morphisms
$F_{N^*(\Delta^{m'}\times\Delta^{n'})}\,\rightarrow\,N^*(\Delta^{m'})\otimes N^*(\Delta^{n'})$
determine a morphism
from the latching algebra $F_{N^*(\Delta^{m'}\times\Delta^{n'})}$
to $N^*(\Delta^m)\otimes N^*(\Delta^n)$.
Precisely,
we consider the composite
$$\colim_{(u,v)} F_{N^*(\Delta^{m'}\times\Delta^{n'})}
\,\rightarrow\,\colim_{(u,v)} N^*(\Delta^{m'})\otimes N^*(\Delta^{n'})
\,\rightarrow\,N^*(\Delta^m)\otimes N^*(\Delta^n).$$
Similarly,
we have a morphism from the cofibrant algebra $F_{N^*(\Delta^m\times\Delta^n)}$
to the matching algebra $\lim_{(u,v)} N^*(\Delta^{m'})\otimes N^*(\Delta^{n'})$.
In this case,
we consider the composite
$$F_{N^*(\Delta^m\times\Delta^n)}
\,\rightarrow\,\lim_{(u,v)} F_{N^*(\Delta^{m'}\times\Delta^{n'})} 
\,\rightarrow\,\lim_{(u,v)} N^*(\Delta^{m'})\otimes N^*(\Delta^{n'}).$$
Finally,
we obtain a commutative diagram
$$\xymatrix{\colim_{(u,v)} F_{N^*(\Delta^{m'}\times\Delta^{n'})}\ar[rr]\ar[d] & &
N^*(\Delta^m)\otimes N^*(\Delta^n)\ar[d] \\
F_{N^*(\Delta^m\times\Delta^n)}\ar[rr]\ar@{-->}[rru] & &
\lim_{(u,v)} N^*(\Delta^{m'})\otimes N^*(\Delta^{n'}) \\ }$$
We deduce the existence of a fill-in morphism from the results of lemma \ref{ReedyAlgebras}.
We claim that this morphism satisfies our requirements.

Precisely,
by construction,
the functoriality diagrams
$$\xymatrix{ F_{N^*(\Delta^{m'}\times\Delta^{n'})}\ar[d]\ar[rr] & &
N^*(\Delta^{m'})\otimes N^*(\Delta^{n'})\ar[d] \\
F_{N^*(\Delta^m\times\Delta^n)}\ar[rr] & &
N^*(\Delta^m)\otimes N^*(\Delta^n) \\ }$$
commute
for pairs of surjective morphisms
$u: [m]\,\rightarrow\,[m']$ and $v: [n]\,\rightarrow\,[n']$.
The general case follows from the induction assumption,
because a morphism in the simplicial category
has a factorization $u = u' u''$,
where $u'$ is injective and $u''$ is surjective.
Similarly,
we prove that our morphism
$F_{N^*(\Delta^m\times\Delta^n)}\,\rightarrow\,N^*(\Delta^m)\otimes N^*(\Delta^n)$
makes all functoriality diagrams commute
for pairs of morphisms
$u: [m']\,\rightarrow\,[m]$ and $v: [n']\,\rightarrow\,[n]$,
such that $m'+n'<N$.

This achieves the proof of lemma \ref{EZModel}.

\tit{\it Construction of the homotopy Eilenberg-Zilber equivalence}\label{EZConstruction}\endtit
We consider the category formed by pairs of morphisms $\sigma: \Delta^m\,\rightarrow\,X$
and $\tau: \Delta^n\,\rightarrow\,Y$.
We have an induced morphism
$(\sigma\times\tau)^*: F_{X\times Y}\,\rightarrow\,F_{\Delta^m\times\Delta^n}$
for each pair $\sigma\times\tau: \Delta^m\times\Delta^n\,\rightarrow\,X\times Y$.
We consider the composites of these morphisms
with the natural transformation
$F_{\Delta^m\times\Delta^n}\,\rightarrow\,N^*(\Delta^m)\otimes N^*(\Delta^n)$
supplied by lemma \ref{EZModel}.
Since our constructions are all functorial,
we obtain a natural morphism
$$F_{X\times Y}\,\rto{}{}\,\lim_{(\sigma,\tau)} N^*(\Delta^m)\otimes N^*(\Delta^n).$$
We recall that $X = \colim_{\sigma: \Delta^m\,\rightarrow\,X} \Delta^m$
and $Y = \colim_{\tau: \Delta^n\,\rightarrow\,Y} \Delta^n$
({\it cf.} [\bibref{GoJ}]).
Consequently,
we have
$$N^*(X)\widehat{\otimes} N^*(Y) = \lim_{(\sigma,\tau)} N^*(\Delta^m)\otimes N^*(\Delta^n)$$
and the construction above provides a morphism
$F_{X\times Y}\,\rightarrow\,N^*(X)\widehat{\otimes} N^*(Y)$,
which is functorial in $X$ and $Y$.
For $X = \Delta^m$ and $Y = \Delta^n$,
this morphism agrees clearly with the map
$F_{\Delta^m\times\Delta^n}\,\rightarrow\,N^*(\Delta^m)\otimes N^*(\Delta^n)$
defined by lemma \ref{EZModel}.
In particular, for $X = Y = \Delta^0$,
we can assume that our map $F_{\Delta^0\times\Delta^0}\,\rightarrow\,N^*(\Delta^0)\otimes N^*(\Delta^0)$
is given by the augmentation morphism $F_{\F}\,\rightarrow\,\F$.

The proof of lemma \ref{EZDefinition} is complete.

\references

\medskip{\parindent=0.5cm\leftskip=0cm

\biblabel{BE}\refto{M. Barratt, P. Eccles},
{\it On $\Gamma_+$-structures. I. A free group functor for stable homotopy theory},
Topology {\bf 13} (1974), 25-45.

\biblabel{BF}\refto{C. Berger, B. Fresse},
{\it Combinatorial operad actions on cochains},
preprint {\tt arXiv: math.AT/0109158} (2001).

\biblabel{FS}\refto{E. Dror-Farjoun, J. Smith},
{\it A geometric interpretation of Lannes' functor $T$},
in {``International Conference on Homotopy Theory (Marseille-Luminy, 1988)''},
Ast\'erisque {\bf 191} (1990), 87-95.

\biblabel{D}\refto{W. Dwyer}, {\it Strong convergence of the Eilenberg-Moore spectral sequence},
Topology {\bf 13} (1974), 255-265.

\biblabel{DEZ}\reftosame, {\it Homotopy operations for simplicial commutative algebras},
Trans. Amer. Math. Soc. {\bf 260} (1980), 421-435. 

\biblabel{DHK}\refto{W. Dwyer, P. Hirschhorn, D. Kan},
{\it Model categories and general abstract homotopy theory},
in preparation.

\biblabel{F}\refto{B. Fresse}, {\it Lie theory of formal groups over an operad},
J. Algebra {\bf 202} (1998), 455-511.

\biblabel{GeJ}\refto{E. Getzler, J. Jones},
{\it Operads, homotopy algebra and iterated integrals for double loop spaces},
preprint {\tt arXiv:hep-th/9403055} (1994).

\biblabel{GiK}\refto{V. Ginzburg, M. Kapranov},
{\it Koszul duality for operads}, Duke Math. J. {\bf 76} (1994), 203-272.

\biblabel{GoJ}\refto{P. Goerss, J. Jardine},
Simplicial homotopy theory, Progress in Mathematics {\bf 174}, Birkh\"auser, 1999.

\biblabel{H}\refto{M. Hovey}, Model categories,
Mathematical Surveys and Monographs {\bf 63}, American Mathematical Society, 1998.

\biblabel{L}\refto{J. Lannes},
{\it Sur les espaces fonctionnels dont la source est le classifiant d'un $p$-groupe ab\'elien \'el\'ementaire},
Publ. Math. Inst. Hautes Etudes Sci. {\bf 75} (1992), 135-244. 

\biblabel{MLH}\refto{S. Mac Lane}, Homology, Die Grundlehren der mathematischen Wissenschaften {\bf 114},
Springer-Verlag, 1963.

\biblabel{MLCat}\reftosame, Categories for the working mathematician, Second Edition,
Graduate Texts in Mathematics {\bf 5}, Springer Verlag, 1998.

\biblabel{Mnd}\refto{M. Mandell},
{\it $E_\infty$ algebras and $p$-adic homotopy theory}, Topology {\bf 40} (2001), 43-94.

\biblabel{Mrl}\refto{F. Morel},
{\it Ensembles profinis simpliciaux et interpr\'etation g\'eom\'etrique du foncteur $T$},
Bull. Soc. Math. France {\bf 124} (1996), 347-373. 

\biblabel{Q}\refto{D. Quillen}, Homotopical algebra, Lecture Notes in Mathematics {\bf 43}, Springer-Verlag, 1967.

\biblabel{R}\refto{C.L. Reedy}, {\it Homotopy theory of model categories}, unpublished manuscript (1973).

\biblabel{Sch}\refto{L. Schwartz},
Unstable modules over the Steenrod algebra and Sullivan's fixed point set conjecture,
Chicago Lectures in Mathematics, University of Chicago Press, 1994.

\biblabel{Stn}\refto{N. Steenrod},
{\it Products of cocycles and extensions of mappings},
Ann. of Math. {\bf 48} (1947), 290-320.

\biblabel{SmrCLoop}\refto{V. Smirnov}, {\it On the chain complex of an iterated loop space},
Izv. Akad. Nauk SSSR Ser. Mat. {\bf 53} (1989), 1108-1119.
English translation in Math. USSR-Izv. {\bf 35} (1990), 445-455.

\biblabel{SmrHLoop}\reftosame, {\it The homology of iterated loop spaces},
preprint {\tt arXiv:math.AT/0010061} (2000).

}

\medskip{\it Acknowledgements:}
I have learned about Mandell's work from a course given by Paul Goerss
and
I am grateful to Ezra Getzler and Paul Goerss for their invitation
at Northwestern University in spring 1999.

\bye